.
\font\piccolo=cmr8.
\font\piccolissimo=cmr6.
\font\script=eusm10.
\font\sets=msbm10.
\font\stampatello=cmcsc10.
.
\def\0{{\bf 0}}
\def\1{{\bf 1}}
\def\defineq{\buildrel{def}\over{=}}
\def\EqByDef{\buildrel{\bullet}\over{=}}
\def\definiz{\buildrel{def}\over{\Longleftrightarrow}}
\def\C{\hbox{\sets C}}

\def\N{\hbox{\sets N}}
\def\R{\hbox{\sets R}}
\def\Primes{\hbox{\sets P}}

\def\Z{\hbox{\sets Z}}
\def\Ide{\mathop{\hbox{\rm Id}}}
\def\Pb{(P)^{\flat}}
\def\square{\hbox{\vrule\vbox{\hrule\phantom{s}\hrule}\vrule}}
\def\nondivide{\not |\thinspace }
\def\supporto{{\rm supp}}
\def\AFwin{\C^{\N}_{{\rm Win}}}
\def\AFsmwin{\C^{\N}_{{\rm smWin}}}
\def\AFfinwin{\C^{\N}_{{\rm finWin}}}
\def\WOD{(\hbox{\stampatello WOD})}
\def\WAD{(\hbox{\stampatello WAD})}
\def\WA{(\hbox{\stampatello WA})}
\def\DH{(\hbox{\stampatello DH})}

\def\Reef{\hbox{\stampatello REEF}}
\def\WReef{\hbox{\stampatello Weak REEF}}

\def\Irr{{\hbox{\rm Irr}}}
\def\IrrPdF{{\hbox{\rm Irr}^{(P)}_d\,F}}
\def\IrrP1F{{\hbox{\rm Irr}^{(P)}_1\,F}}
\def\IrrPF{{\hbox{\rm Irr}^{(P)}\,F}}

\def\Car{{\rm Car}}
\def\Win{{\rm Win}}
\def\CarqF{{\rm Car}_q\,F}
\def\WinqF{{\rm Win}_q\,F}
\def\CarT{\Car\; }
\def\WinT{\Win\; }
\def\CarTP{\Car^{(P)}\, }
\def\WinTP{\Win^{(P)}\, }
\def\CarPqF{\Car^{(P)}_q F\;}
\def\WinPqF{\Win^{(P)}_q F\;}

\def\WinPbqF{\Win^{\Pb}_q F\;}
\def\Rvl{(\hbox{\stampatello Rvl})}
\def\NSL{(\hbox{\stampatello NSL})}
\def\IPP{(\hbox{\stampatello IPP})}
\def\ETD{(\hbox{\stampatello ETD})}

\def\QED{\hfill {\rm QED}\par}

\def\assurdo{\rightarrow\!\leftarrow}
\def\sumflat{\mathop{{\sum}^{\flat}}}

\def\RamaSer{{\mathop{\hbox{\script R}}}}
\def\RamaSmSer{{\mathop{\widetilde{\hbox{\script R}}}}}
\par
\centerline{\bf On Ramanujan smooth expansions for a general arithmetic function}
\bigskip
\centerline{Giovanni Coppola}\footnote{ }{MSC $2020$: $11{\rm N}37$ - Keywords: Ramanujan expansion, smooth numbers, arithmetic functions} 

\bigskip

\par
\noindent
{\bf Abstract}. We study in detail the Ramanujan smooth expansions, for arithmetic functions; we start with the most general ones, for which we supply the \lq \lq $P-$local expansions\rq \rq, for arguments with all prime-factors $p\le P$ (namely, $P-$smooth arguments), that are also square-free; then, we supply general results for interesting subsets of arithmetic functions, regarding both their $P-$local and (global) Ramanujan smooth expansions. 

\bigskip

\par
\noindent{\bf 1. Introduction. Definitions and Notations. A local expansion for all arithmetic functions} 
\bigskip
\par
\noindent
The history-making paper by Ramanujan [R] studies the {\stampatello Ramanujan Sum} of {\it modulus} $q$ and {\it argument} $a$ :
$$
\forall q\in \N,\enspace \forall a\in \Z,
\qquad
c_q(a)\defineq \sum_{{j\le q}\atop {(j,q)=1}}\cos {{2\pi ja}\over q} = \sum_{j\in \Z_q^*}e_q(ja), 
$$
\par
\noindent
where $(j,q)$ is the greatest common divisor of $j$ and $q$,\enspace $\Z_q^*$\enspace is the set of {\it reduced residue classes modulo} $q$\enspace and 
\enspace $\forall q\in \N,$ \thinspace $\forall m\in \Z,$ \enspace $e_q(m)\defineq e^{2\pi im/q}.$ In passing, see that $c_q(0)=\varphi(q)\defineq |\Z_q^*|$ is {\it Euler's Totient}. 
\quad These sums were already known before, but Ramanujan used them in [R] in His \lq \lq {\stampatello Ramanujan expansions}\rq \rq, for some arithmetic functions. Of course, He didn't study these expansions in full theoretical detail, but giving an illuminating example: He expanded the {\it null-function} $\0$, i.e., $\0(n)\defineq 0$, $\forall n\in \N$, as
$$
\0(a)=\sum_{q=1}^{\infty}{1\over q}c_q(a)\defineq \lim_x \sum_{q\le x}{1\over q}c_q(a), 
\qquad 
\forall a\in \N. 
$$
\par
\noindent
These \lq \lq {\stampatello Ramanujan coefficients}\rq \rq, say $R_0(q)\defineq 1/q$, $\forall q\in \N$, for $\0$, {\stampatello are not unique}, as we'll see soon. 
\par
We call, once the arithmetic function $G:\N \rightarrow \C$ is fixed, 
$$
\forall a\in \N,
\quad
\RamaSer_G(a)\defineq \lim_x \sum_{q\le x}G(q)c_q(a),
\enspace
\hbox{\stampatello the\enspace Ramanujan\enspace Series\enspace of\enspace coefficient} \enspace G,
$$
\par
\noindent
whenever, $\forall a\in \N$, this limit exists in $\C$, i.e., the series converges pointwise in all $a\in \N$. 
\par
Notice the \lq \lq {\stampatello classic summation}\rq \rq: $q\le x$. \enspace (In the following, $\N_0\defineq \N \cup \{0\}$ and $\Primes$ is the set of primes.) 
\par
\noindent
An immediate consequence is : 
$$
\RamaSer_G=\RamaSer_{G+\alpha R_0},
\qquad 
\forall \alpha \in \C,
$$
\par
\noindent
in other words, if the arithmetic function $F$ has the Ramanujan expansion $F=\RamaSer_G$, for some $G:\N \rightarrow \C$, then $F$ has, also, {\stampatello an infinity} (as $\alpha \in \C$ varies) {\stampatello of other Ramanujan coefficients}. See that all of this is implicit in Ramanujan's paper [R]. Actually, Hardy [H] gave the independent $H_0(q)\defineq 1/\varphi(q)$ for $\0=\RamaSer_{H_0}$; however, we have many-more than $\alpha H_0 +\beta R_0$ as Ramanujan coefficients for $\0$, beyond $R_0$ : [CG1].
\par
Ramanujan considered expansions with the above, say, {\stampatello classic} method of summation. 
\par
\noindent
We introduce a {\stampatello smooth} summation method: our partial sums are not over $q\le x$, for fixed $x$, but over $q\in (P)$, for fixed $P\in \Primes$, a prime, where $(P)\defineq \{n\in \N : p|n \Rightarrow p\le P\}$ is the set of $P-${\it smooth} numbers: 
$$
\forall a\in \N 
\enspace \hbox{\stampatello and} \enspace 
\forall P\in \Primes
\enspace \hbox{\stampatello fixed}, \enspace 
\enspace \hbox{\stampatello the $P-$smooth partial sum} \enspace 
\sum_{q\in (P)}G(q)c_q(a)
\enspace \hbox{\stampatello is finite,} \enspace 
\forall G:\N \rightarrow \C.
$$
\par
\noindent
This follows from the property called \lq \lq {\stampatello Ramanujan vertical limit}\rq \rq: with $v_p(q)\defineq \max\{K\in \N_0 : p^K|q \}$ the $p-$adic {\it valuation} of $q\in \N$, $\forall p\in \Primes$, 
$$
c_q(a)\neq 0
\enspace \Longrightarrow \enspace
v_p(q)\le v_p(a)+1,
\enspace
\forall p|q.
$$
\par				
\noindent
For a Proof of a slightly stronger version (with the equivalence!), see in $\S3.1$. 

\medskip

\par
\noindent
We introduce, once the arithmetic function $G:\N \rightarrow \C$ is fixed, 
$$
\forall a\in \N,
\quad
\RamaSmSer_G(a)\defineq \lim_P \sum_{q\in (P)}G(q)c_q(a),
\enspace
\hbox{\stampatello the\enspace Ramanujan\enspace Smooth\enspace Series\enspace of\enspace coefficient} \enspace G,
$$
\par
\noindent
whenever, $\forall a\in \N$, this limit over primes $P\to \infty$ exists in $\C$, i.e., $P-$smooth partial sums converge pointwise in all $a\in \N$. 
\par
See that, in general, {\stampatello classic} Ramanujan series are different from {\stampatello smooth} Ramanujan series and this can't be appreciated, of course, when we are {\it in} the {\it case of absolute convergence}: 
$$
\sum_{q=1}^{\infty}\left| G(q)c_q(a)\right|<\infty. 
$$ 
\par
\noindent
In fact, as we know, given a fixed $a\in \N$, this property of {\bf absolute convergence of Ramanujan series} in the argument $a$ {\stampatello entails the same sum, for all summation methods}, of our series in $a$.
\par
However, when the absolute convergence of Ramanujan series is not guaranteed, we may have not only different values for $\RamaSer_G(a)$ and $\RamaSmSer_G(a)$, but even that at least one of them doesn't exist : see following example. 
\par
We saw above that $R_0(q)=1/q$ was given by Ramanujan himself as a coefficient, for the null-function. This, of course, with classic summation: $\0=\RamaSer_{R_0}$.
\par
Well, we prove now the example of Ramanjuan smooth coefficient $\1$, the constant-$1$-function: $\1(n)\defineq 1$, $\forall n\in \N$, for our arithmetic function $\0$; i.e., our first {\stampatello Ramanujan smooth expansion}, namely $\0=\RamaSmSer_{\1}$, is: 
$$
\0(a)=\RamaSmSer_{\1}(a)
 =\lim_P \sum_{q\in (P)}\1(q)c_q(a)
  =\lim_P \sum_{q\in (P)}c_q(a)
   =\lim_P \prod_{p\le P}\sum_{K=0}^{\infty}c_{p^K}(a),
    =\lim_P \prod_{p\le P}\sum_{K=0}^{v_p(a)+1}c_{p^K}(a),
\enspace 
\forall a\in \N, 
$$
\par
\noindent
because we use: fixed $a$, $c_q(a)$ is multiplicative w.r.t.(with respect to) $q$ and, then, we apply Ramanujan vertical limit we saw above. Now: ALL these $P-$smooth partial sums vanish, as we see from H\"older's Formula [M] with $\mu$ the {\it M\"obius function}, compare [C5] basic facts,
$$
c_q(a)=\varphi(q){{\mu(q/(q,a))}\over {\varphi(q/(q,a))}}
\enspace \Longrightarrow \enspace
c_{p^K}(a)=c_{p^K}(p^{v_p(a)})=\varphi(p^K){{\mu(p^{K-\min(K,v_p(a))})}\over {\varphi(p^{K-\min(K,v_p(a))})}}
$$
\par
\noindent
entailing 
$$
c_{p^K}(a)=\varphi(p^K),
\thinspace \forall 0\le K\le v_p(a), \enspace
c_{p^{v_p(a)+1}}(a)=-p^{v_p(a)}
\enspace \Longrightarrow \enspace
\sum_{K=0}^{v_p(a)+1}c_{p^K}(a)=\sum_{K=0}^{v_p(a)}\varphi(p^K)-p^{v_p(a)}=0
$$
\par
\noindent
and valid for all fixed $a\in \N$ (once again, convergence is pointwise in $a\in \N$) : {\stampatello finite Euler product} above, over $p\le P$, has ALL $p-$factors vanishing and we are done. 
\par
Of course, this series doesn't converge absolutely: the partial sums
$$
\sum_{q\in (P)}|c_q(a)|=\prod_{p\le P}\sum_{K=0}^{v_p(a)+1}|c_{p^K}(a)|
 =\prod_{p\le P}\left(\sum_{K=0}^{v_p(a)}\varphi(p^K)+p^{v_p(a)}\right)
  =\prod_{p\le P}(2p^{v_p(a)})
   \ge \prod_{p\le P}2
$$
\par
\noindent
are bounded from below by $2^{\pi(P)}$, with $\pi(x)\defineq |\{p\in \Primes : p\le x\}|$, so the limit over $P\to \infty$ is $+\infty$. 
\par
While $\RamaSmSer_{\1}=\0$, the classic $\RamaSer_{\1}$ DOESN'T EXIST, as the general term $c_q(a)$ is NOT INFINITESIMAL as $q\to \infty$ (because, $\forall a\in \N$ fixed, $c_q(a)\in \Z$, as we know from H\"older's formula above, so that $c_q(a)\to 0$, as $q\to \infty$, means: $\exists Q\in \N$, $c_q(a)=0$, $\forall q>Q$ and THIS goes against Ramanujan sums vanishing, see $\S3.1$). 

\medskip

\par
Most weird things are also most interesting things, so to speak: they might happen ONLY in absence of absolute convergence! 

\vfill
\eject

\par				
\noindent
We leave aside, for the moment, convergence issues and we focus for a while on properties and definitions regarding arithmetic aspects. 
\par
First of all, the most fundamental property we apply, regarding Ramanujan sums, when coming to any kind of expansion in this paper, is {\stampatello Ramanujan Orthogonality} : 
$$
\1_{q|a}=\1_{a\equiv 0(\!\!\bmod q)}={1\over q}\sum_{\ell | q}c_{\ell}(a), 
\leqno{(1.1)}
$$
\par
\noindent
valid for all \lq \lq moduli\rq \rq, say, $q\in \N$ and integer arguments $a$; hereafter we abbreviate with $\1_{\wp}$ the indicator function of property $\wp$ : it's $\defineq 1$ IFF (If and only If) $\wp$ is true and $\defineq 0$ IFF $\wp$ is false; actually, we also use $\1_{S}$, for any subset of natural numbers $S\subseteq \N$, as the indicator function of $S$: we saw above $\1$, namely $\1=\1_{\N}$. The name we use here for above property is, of course, in honor of Ramanujan; however, there will be other kinds of Orthogonality, for Ramanujan sums, along this paper : compare well-known Facts in [C5]. 
\par
By the way, we take for granted ALL the {\bf basic facts} (like the Orthogonality of Additive Characters) quoted in [C5] and quotations therein of Bibliographic sources (i.e., Davenport's [D] \& Tenenbaum's [T] Books, with Iwaniec-Kowalski [IKo] Book, too, see [C5]). 
\par
Here, compare [IKo] notation, we write $\sumflat_{\ldots}$ for the sum $\sum_{\ldots}$ with restriction to {\it square-free} summands (with other constraints denoted $\ldots$ here). As a kind of extension for this notation, we add the {\it superscript} $\flat$ to any subset of natural numbers, say $S$, writing $S^{\flat}$ for the square-free elements of $S$; in particular, $\N^{\flat}$ is the subset of square-free natural numbers (recall: $\N\defineq \{1,2,3,4,\ldots\}$), like $\Pb\defineq \{ n\in (P) : n\enspace $ square-free $\}$ is the subset of $P-$smooth numbers that are square-free, too. Notice that, once fixed $P\in \Primes$, $\Pb$ is a finite set, of cardinality $|\Pb|=2^{\pi(P)}$, with $\max\Pb=\prod_{p\le P}p$ the $P-${\bf primorial}. 
\par 
We add, to our notation, the set of $P-${\it sifted naturals}: 
$$
)P(\defineq \{ n\in \N : p|n\enspace \Rightarrow \enspace p>P\}
\EqByDef \{ n\in \N : (n,p)=1, \forall p\le P\}
\EqByDef \{ n\in \N : (n,\prod_{p\le P}p)=1\}, 
$$
\par
\noindent
hereafter indicating with $\EqByDef$ that there is equality by definition. Notice: last way we write it highlights the sieving procedure, to identify its elements (compare the elementary Eratosthenes-Legendre Sieve [HaRi]). 
\par
Once fixed $P\in \Primes$, not only $(P)\cap )P(=\{1\}$, but an important feature is that, even if $(P)\cup )P($ is NOT all $\N$, we have: 
$$
\forall a\in \N, a=rs,
\quad
s\in (P),\enspace r\in )P(, 
$$
\par
\noindent
for {\bf unique factors} $a_{(P)}\defineq s\EqByDef\prod_{p\le P}p^{v_p(a)}$ (for: $s$mooth) and $a_{)P(}\defineq r\EqByDef\prod_{p>P}p^{v_p(a)}$ (for: $r$ough); this kind of decomposition was used at the state-of-art by Paul Erd\"os in His Proofs. Here, much more modestly, we'll see this kind of {\bf orthogonality} (since $r$ \& $s$ are coprime, a kind of ortogonality) at work, so to speak, in our next Theorem 2.1, named in honor of Wintner: Wintner's Orthogonal Decomposition, say, $\WOD$. (Compare $\S3.1.4$.)

\bigskip

\par
We give now an elementary result: a $P-${\it local expansion}, for {\bf all} arithmetic functions. Recall $F'$ is the {\stampatello Eratosthenes Transform} of our arithmetic function $F$ (see [T], compare [C5]).
\smallskip
\par
\noindent{\bf Theorem 1.1}. {\it Let } $F:\N \rightarrow \C$. {\it Fix a prime } $P$. {\it Then, defining } ${\displaystyle \WinPbqF \defineq \sum_{{d\in \Pb}\atop {d\equiv 0(\!\! \bmod q)}}{{F'(d)}\over d} }$, $\forall q\in \N$, {\it we have the $P-$local expansion}  
$$
F(a)=\sum_{q\in \Pb}\left(\WinPbqF\right)c_q(a),
\quad 
\forall a\in \Pb. 
$$

\vfill

\par
\noindent
This is a kind of coming soon, for Lemma $3.2$, in which we give many more details and the one-line-proof: it uses ONLY $(1.1)$ above. By the way, in THIS paper we apply \& quote Lemma $3.2$, while OUTSIDE present paper we'll quote Theorem 1.1! 
\par
It's a kind of UNIVERSAL expansion, namely it holds FOR ALL ARITHMETIC FUNCTIONS. Of course, (1): it's local, namely both the coefficients and their indices (outside which they vanish, see $\S3$) depend on $P$; (2): it's only for square-free and $P-$smooth arguments; (3):last, but not least, we have no idea of what are these coefficients in this $P-$local expansion. By the way, the definition above for the $d-$series is well-posed as: $d\in \Pb$ is in a finite set, so the $d-$series converges trivially, being a finite sum!  

\eject

\par				
\noindent{\bf 2. Main Results: Statements and Proofs} 
\bigskip
\par
\noindent
We collect here all the Statements \& Proofs of our main results, but the ones needing the Lemmata in $\S3$ \& $\S4$ will be proved in next subsection $\S2.1$. 
\medskip
\par
By the way, we call \lq \lq main\rq \rq, the present results, since also the Facts, Lemmata and Propositions given in $\S3$ and $\S4$, $\S5$ have their own interest. 
\medskip
\par
We state and prove our next result. Compare [HaRi] for the Eratosthenes-Legendre Sieve. 
\smallskip
\par
\noindent{\bf Theorem 2.1}.
%
%
 {\stampatello (Wintner Orthogonal Decomposition for $F'$)} 
\par
\noindent
{\it Let } $F:\N \rightarrow \C$ {\it have the} {\stampatello Wintner Transform}, {\it say}
$$
\WinT F : q\in \N \mapsto \WinqF \defineq \sum_{d\equiv 0(\!\!\bmod q)}{{F'(d)}\over d}\in \C, 
$$
{\it i.e., has all such} {\stampatello Wintner's $q-$th coefficients}: {\it these series converge pointwise } $\forall q\in \N$. {\it Fix a prime} $P$. {\it Then}
$$
\forall d\in \N, 
\qquad
\IrrPdF \defineq \sum_{{r\in )P(}\atop {r>1}}{{F'(dr)}\over r}
 \EqByDef \lim_x \sum_{{r\in )P(}\atop {1<r\le x}}{{F'(dr)}\over r}
$$
\par
\noindent
{\it is a series, converging pointwise in } $d$, $\forall d\in \N$, {\it called the } {\stampatello Irregular Series} {\it of our } $F$, {\stampatello with argument } $d$, {\stampatello over the prime } $P$. {\it Varying $d\in \N$, we get as arithmetic function the } $P-${\stampatello Irregular Series} {\it of } $F$, 
$$
\IrrPF : d\in \N \mapsto \sum_{{r\in )P(}\atop {r>1}}{{F'(dr)}\over r}\in \C. 
$$
\par
\noindent
{\it In all, } $\exists \WinT F$ {\it entails } $\exists \IrrPF$, $\forall P\in \Primes$; {\it and the} {\stampatello Decomposition of Eratosthenes Transform} $F'$ : $\forall P\in \Primes$,
$$
F'(d)=d\sum_{K\in (P)}\mu(K)\Win_{dK}\,F - \IrrPdF,
\qquad
\forall d\in \N. 
\leqno{\WOD_{F',P}}
$$
\par
\noindent
{\it Thus in particular}
$$
\forall d\in \N, 
\qquad
F'(d)=\lim_P\left(d\sum_{K\in (P)}\mu(K)\Win_{dK}\,F - \IrrPdF\right). 
\leqno{\WOD_{F'}}
$$
\par
\noindent{\bf Proof}. We prove both the {\stampatello convergence} of $\IrrPdF$ and $\WOD_{F',P}$ : {\stampatello fix} $P\in \Primes$ and $d\in \N$. Then, $r\in )P($ is equivalent to: $(r,\prod_{p\le P}p)=1$ (see the $P-$sifted set's definition on page 3), which can be expressed as: ${\displaystyle \1_{r\in )P(}=\sum_{K|r,K|\prod_{p\le P}p}\mu(K)=\sum_{{K|r}\atop {K\in(P)}}\mu(K) }$, compare [C5], entailing, for a large $x>0$, 
$$
\sum_{{r\in )P(}\atop {r\le x}}{{F'(dr)}\over r}=\sum_{r\le x}\sum_{{K|r}\atop {K\in(P)}}\mu(K){{F'(dr)}\over r}
 =\sum_{K\in(P)}\mu(K)\sum_{m\le x/K}{{F'(dKm)}\over {Km}}
  =d\sum_{K\in(P)}\mu(K)\sum_{m\le x/K}{{F'(dKm)}\over {dKm}}, 
$$
\par
\noindent
whence, {\stampatello exchanging} the {\stampatello finite sum over} $K\in \Pb$ {\stampatello with} $\lim_x$, 
$$
F'(d)+\lim_x \sum_{{r\in )P(}\atop {1<r\le x}}{{F'(dr)}\over r}=d\sum_{K\in(P)}\mu(K)\lim_x \sum_{{t\le dx}\atop {t\equiv 0(\!\!\bmod dK)}}{{F'(t)}\over t}
 =d\sum_{K\in(P)}\mu(K)\Win_{dK}\,F, 
$$
\par
\noindent
from the definition of \enspace $\Win_{dK}\,F$, guaranteed by hypothesis: $\exists \WinT\,F$.\hfill $\square$

\medskip 

\par				
\noindent
We state and prove our next result, an immediate consequence of Theorem 2.1. 
\smallskip
\par
\noindent{\bf Corollary 2.1}.
%
%
 {\stampatello (Wintner Orthogonal Decomposition for $F$)} 
\par
\noindent
{\it Let } $F:\N \rightarrow \C$ {\it have the} {\stampatello Wintner Transform}. {\it Fix a prime } $P$. {\it Then}
$$
F(a)=\sum_{q\in (P)}\left(\Win_q\,F\right)c_q(a) - \sum_{d|a}\IrrPdF,
\qquad
\forall a\in (P). 
\leqno{\WOD_{F,P}}
$$
\par
\noindent
{\it Thus in particular}
$$
\forall a\in \N, 
\quad
F(a)=\lim_P\left(\sum_{q\in (P)}\left(\Win_q\,F\right)c_q(a) - \sum_{d|a}\IrrPdF\right). 
\leqno{\WOD_{F}}
$$
\par
\noindent{\bf Proof}. The rough idea is to sum $\WOD_{F',P}$ over divisors $d|a$. Right approach follows. We apply Kluyver's formula [K] (see [C5] basic facts), recall: 
$$
c_q(a)=\sum_{{d|a}\atop {d|q}}d\mu\left({q\over d}\right), 
$$ 
\par
\noindent
to above decomposition for $F'$, i.e. $\WOD_{F',P}$ , for a fixed $P\in \Primes$, getting $\WOD_{F,P}$ , because we set $K=q/d$ and we use the elementary fact: $a\in (P)$ and $d|a$ entail $d\in (P)$; namely, $P-$smooth numbers are a {\bf divisor-closed} set.\hfill $\square$

\medskip

\par
Wintner Orthogonal Decomposition is very simple and beautiful; it starts from the most elementary sieve method, namely the inclusion-exclusion principle with products of primes, codified by the M\"obius function, [HaRi] : it is, in fact, called the Sieve of Eratosthenes-Legendre and, not by chance, we are working with the Eratosthenes Transform, here (name due to Wintner in his beautiful Book [W], by the way see Appendix $\S5$ here). 
\par
It may seem incredible, but it turns very useful for our purposes: we are seeking Ramanujan Smooth Expansions. And $\WOD_F$ in previous Corollary not only gives a beautiful, rather unbelievable characterization of such Expansions in terms of, say, a series summed over sifted numbers, that's our Irregular Series for $F$; also, as a gift, these Smooth Expansions have a completely explicit coefficient: Wintner's Transform of our $F$. All of this sparkles, from Corollary 2.1, in our next Corollary 2.2!   

\medskip

\par
We state and prove our next result, an immediate consequence of Corollary 2.1. 
\smallskip
\par
\noindent{\bf Corollary 2.2}.
%
%
 {\stampatello (Characterizing Ramanujan Smooth Expansion with Wintner coefficients)} 
\par
\noindent
{\it Let } $F:\N \rightarrow \C$ {\it have } $\WinT\,F$. {\it Then}
$$
\forall a\in \N,
\enspace
F(a)=\lim_P \sum_{q\in (P)}\left(\Win_q\,F\right)c_q(a)
\quad
\Longleftrightarrow
\quad
\forall d\in \N,
\enspace
\lim_P \IrrPdF =0
\quad
\Longleftrightarrow
\quad
\lim_P \IrrPF = \0. 
$$
\par
\noindent{\bf Proof}. Of course, second equivalence holds as a definition, so we apply $\WOD_{F}$ to prove first equivalence: 
$$
\forall a\in \N,
\thinspace
F(a)=\lim_P \sum_{q\in (P)}\left(\Win_q\,F\right)c_q(a)
\enspace 
\Longleftrightarrow
\enspace 
\forall a\in \N,
\thinspace
\sum_{d|a} \lim_P \IrrPdF = 0
\enspace 
\Longleftrightarrow
\enspace 
\forall d\in \N,
\thinspace
\lim_P \IrrPdF =0. 
$$
\par
\noindent
After the above exchange of $\lim_P$ with the finite sum $\sum_{d|a}$ (independent of $P$, too), we simply have applied M\"obius Inversion, [T], compare [C5].\hfill $\square$

\medskip

\par
Since $\C^{\N}$ indicates ALL arithmetic functions, arithmetic functions $F$ having the Wintner Transform will be indicated as: $\AFwin$; furthermore\enspace $\AFsmwin \defineq \{ F\in \AFwin : \exists P\in \Primes, \supporto(\WinT\,F)\subseteq (P)\}$ are the $F$ with $\WinT\,F$, say, \lq \lq smooth-supported\rq \rq; and\enspace $\AFfinwin\defineq \{ F\in \AFwin : |\supporto(\WinT\,F)|<\infty \}$ are the $F$ with $\WinT\,F$, say, \lq \lq finitely-supported\rq \rq: namely, these $F$ have only a finite number of non-zero Wintner coefficients, eventually none whenever $\WinT\,F=\0$. 
\vfill
\par
We'll be back on this, after next Corollary. 
\eject

\par				
\noindent
In order to expose next result, we need to recall {\stampatello the $q-$th Carmichael coefficient} of our $F:\N \rightarrow \C$, 
$$
\Car_q\,F \defineq {1\over {\varphi(q)}}\lim_x {1\over x}\sum_{a\le x}F(a)c_q(a),
\qquad
\forall q\in \N, 
$$
\par
\noindent
of course whenever, fixed $q\in \N$, this limit exists in the complex numbers. Notice that we defined the Wintner $q-$th coefficient inside Theorem 2.1  above and we also gave there the definition of Wintner Transform of our $F$, whenever ALL $\Win_q\,F$ exist. In the same way, assuming ALL $\Car_q\,F$ exist, we define: 
$$
\CarT\,F : q\in \N\mapsto \Car_q\,F\in \C, 
$$
\par
\noindent
to be {\stampatello the Carmichael Transform} of our $F$. For both these transforms, see $\S3$, in particular $\S3.3$ and compare [C5]. Here we only recall that, under $\WA$, see next Corollary 2.3, these two transforms, not only exist both, but they are also known to be equal: $\CarT\,F=\WinT\,F$. Compare Wintner's Criterion in [ScSp]. 

\medskip

\par
We state and prove our next result. 
\smallskip
\par
\noindent{\bf Corollary 2.3}.
%
%
 {\stampatello (Wintner's \lq \lq Dream\rq \rq\enspace Theorem)} 
\par
\noindent
{\it Let } $F:\N \rightarrow \C$ {\it satisfy} {\stampatello Wintner Assumption}, {\it say}
$$
\sum_{d=1}^{\infty}{{|F'(d)|}\over d}<\infty. 
\leqno{\WA}
$$
{\it Then $F$ has the } {\stampatello Ramanujan Smooth Expansion}, {\it with, say, } {\stampatello Carmichael-Wintner coefficients}
$$
F(a)=\lim_P \sum_{q\in (P)}\left(\CarqF\right) c_q(a)
     =\lim_P \sum_{q\in (P)}\left(\WinqF\right) c_q(a),
\quad
\forall a\in \N. 
\leqno{(2.1)}
$$  
\par
\noindent{\bf Proof}. Thanks to Corollary 2.2, we only need to prove: 
$$
\WA 
\quad \Longrightarrow \quad
\forall d\in \N,
\enspace
\lim_P \IrrPdF =0. 
$$
\par
\noindent
We fix $d\in \N$ and $P\in \Primes$ with $P$ large enough, namely $P>d$. Then, $r\in )P( \, \backslash \, \{1\}$ entails $dr>P$, whence
$$
\IrrPdF \EqByDef d\lim_x \sum_{{r\in )P(}\atop {1<r\le x}}{{F'(dr)}\over {dr}}
\quad \Longrightarrow \quad
\left|\IrrPdF\right|\le d\sum_{m>P}{{|F'(m)|}\over m}
\buildrel{P}\over{\rightarrow} 0,
$$
\par
\noindent
since this $m-$series represents the $P-$tails of the series converging in $\WA$.\hfill $\square$

\medskip

\par
We have now very clear both the importance and the difficulty when handling the Irregular series: the two variables $P\in \Primes$ (say, compare [HaRi], the sifting limit) and $d\in \N$ (our divisors for $F$, the arguments of our $F'$) interact in a non-trivial way! I don't think it would be a good idea, trying to study the uniform convergence of our Irregular Series, w.r.t. the divisors $d\in \N$ ! And likewise it's not a good approach (even if we followed it, in previous Proof!), that of studying the absolute convergence of $\IrrPF$. We are \lq \lq lucky\rq \rq, so to speak, because the strength of Wintner's hypothesis $\WA$ allows not only to get an absolutely convergent $\IrrPdF$, once fixed $d$ and $P$ as above, but also, once fixed $d$,  infinitesimal with $P$ (when I found above Proof, I could not realize its strength \& beauty at once: I still didn't know enough $\IrrPF$ subtleties!).
\smallskip
\par
By the way, previous Proof doesn't allow to prove the absolute convergence of our Smooth Ramanujan Expansion and the reason is clear: already Wintner, in His Book [W], realized that His $\WA$, above, is not sufficient to get a {\stampatello classic Ramanujan Expansion}; this, instead, follows, {\bf with absolute convergence} (so, for any summation method), {\bf from Delange Hypothesis}, recall: $\omega(d)\defineq |\{ p\in \Primes : p|d\}|$, hereafter, 
$$
\sum_{d=1}^{\infty}{{2^{\omega(d)}|F'(d)|}\over d}<\infty. 
\leqno{\DH}
$$
\par
\noindent
The Proof of this is the content of, say, {\it Delange's Theorem} [De] (compare [C1] exposition).
\par
Thus
\smallskip
\par
\noindent{\bf Warning}. Wintner's Dream Theorem ensures ONLY pointwise, not absolute, convergence with the smooth summation method.\hfill $\circ$

\vfill
\eject

\par				
\noindent
We need, for previously defined two subsets of $\AFwin$, two useful definitions; since, of course, $\AFfinwin \subset \AFsmwin$ and this inclusion is strict, we start with $\AFsmwin$, then pass to $\AFfinwin$. 
\par
We define $\forall F\in \AFsmwin$ the {\stampatello Wintner's Prime} for $F$:
$$
P_F \defineq \min\{P\in \Primes : \supporto(\Win\,F)\subseteq (P)\},
\quad
\forall F\in \AFsmwin. 
$$
\par
\noindent
Notice: in particular, for any $F$ with $\WinT\,F=\0$ (so, support's empty) and any $F$ with $\supporto(\WinT\,F)=\{1\}$, we have $P_F=2$. 
\par
While, $\forall F\in \AFfinwin$ with $\WinT\,F=\0$ we set the {\stampatello exact Wintner's Range} for $F$ to: $Q_F\defineq 0$; and it's the following natural number otherwise: 
$$
F\in \AFfinwin, \enspace \WinT\,F\neq \0
\quad \Longrightarrow \quad
Q_F \defineq \max\supporto(\WinT\,F). 
$$

\medskip

\par
We wish to give, in the two environments (one nested in the other), $\AFsmwin$ and $\AFfinwin$, the Wintner Orthogonal Decompositions; these, in fact, for our $F$ in $\AFsmwin$ simplify (and in $\AFfinwin$, even more).
\par
Roughly speaking, $P_F$ is important for $\IrrPF$ and $Q_F$ a must for the sum with Wintner coefficients.

\medskip

\par
We state and prove next Wintner Decompositions for $F'$ \& $F$ into \lq \lq Regular and Irregular Parts\rq \rq, say. 
\smallskip
\par
\noindent{\bf Corollary 2.4}.
%
%
 {\stampatello (Wintner Orthogonal Decompositions when $\WinT\,F$ is smooth-supported)} 
\par
\noindent
{\it Let } $F\in \AFsmwin$, {\it namely Wintner Transform is smooth-supported, with } $P_F$ {\stampatello Wintner's Prime for } $F$. 
{\it Then}
$$
\IrrPF=\Irr^{(P_F)}F,
\quad
\forall P\in \Primes, P\ge P_F. 
\leqno{\hbox{(\stampatello $P-$inertia)}}
$$ 
{\it Thus, previous general formul\ae} $\WOD_{F'}$ {\it and } $\WOD_F$ {\it become } 
$$
\forall F\in \AFsmwin, 
\qquad
F'(d)=d\sum_{K\in (P_F)}\mu(K)\Win_{dK}\,F - \Irr^{(P_F)}_d F, 
\quad
\forall d\in \N
$$
\par
\noindent
{\it and} 
$$
\forall F\in \AFsmwin, 
\qquad
F(a)=\sum_{q\in (P_F)}\left(\WinqF\right)c_q(a) - \sum_{d|a}\Irr^{(P_F)}_d F, 
\quad
\forall a\in \N. 
$$
\par
\noindent
{\it Furthermore, these two formul\ae \enspace specialize, for } $\WinT\,F$ {\stampatello finitely supported}, {\it to}: 
$$
\forall F\in \AFfinwin, 
\qquad
F'(d)=d\sum_{K\le Q_F/d}\mu(K)\Win_{dK}\,F - \Irr^{(P_F)}_d F, 
\quad
\forall d\in \N
$$
\par
\noindent
{\it and} 
$$
\forall F\in \AFfinwin, 
\qquad
F(a)=\sum_{q\le Q_F}\left(\WinqF\right)c_q(a) - \sum_{d|a}\Irr^{(P_F)}_d F, 
\quad
\forall a\in \N. 
$$
\smallskip
\par
\noindent
We leave the Proof to the interested reader.

\medskip

\par
Last but not least, $\hbox{(\stampatello $P-$inertia)}$ above {\stampatello allows to define}, $\forall F\in \AFsmwin$, the 
\par
\leftline{{\stampatello Irregular Part} of our $F$ {\stampatello as} : ${\displaystyle I_F(a)\defineq \sum_{d|a}\Irr^{(P_F)}_d F }$, $\forall a\in \N$, and the
}
\par
\leftline{{\stampatello Smooth Part} of our $F$ {\stampatello as} : ${\displaystyle S_F(a)\defineq \sum_{q\in (P_F)}\left(\WinqF\right)c_q(a) }$, $\forall a\in \N$. 
}
\par
Whenever $F\in \AFfinwin$, we also define the
\par
\leftline{{\stampatello Analytic Part} of our $F$ {\stampatello as} : ${\displaystyle A_F(a)\defineq \sum_{q\le Q_F}\left(\WinqF\right)c_q(a) }$, $\forall a\in \N$. Recall {\bf empty sums are} $\defineq 0$.  
}
\par				
\noindent
Of course, we CAN'T DEFINE $A_F$ for the arithmetic functions $F \in \AFsmwin \backslash \AFfinwin$. For these functions, the Smooth Part has the number of non-zero terms depending on $a\in \N$, by Ramanujan Vertical Limit ! 

\medskip

\par
We give still two definitions, very important for applications (compare [CM], [C5]): for $F\in \AFfinwin$ we say that $F$ {\stampatello has the REEF}, abbreviating {\it Ramanujan Exact Explicit Formula}, exactly when $F=A_F$, I.E.
$$
F(a)=A_F(a)
 \EqByDef \sum_{q\le Q_F}\left(\WinqF\right)c_q(a),
\qquad
\forall a\in \N.
\leqno{(\Reef)}
$$ 
\par
\noindent
When $F\not \in \AFfinwin$, but $F\in \AFsmwin$, we say that $F$ {\stampatello has the Weak REEF} exactly when $F=S_F$, I.E.
$$
F(a)=S_F(a)
 \EqByDef \sum_{q\in (P_F)}\left(\WinqF\right)c_q(a),
\qquad
\forall a\in \N.
\leqno{(\WReef)}
$$ 
\smallskip
\par
\noindent
In both cases, for $F\in \AFfinwin$ to have the $\Reef$, or for $F\in \AFsmwin$ to have the $\WReef$, we need, and it suffices, to have a vanishing Irregular Part for $F$, i.e., $I_F=\0$.
\par
In passing, see that for $F\in \AFfinwin$ we have the regular/irregular decomposition: \enspace $F=A_F-I_F$ \enspace i.e., 
$$
\forall F\in \AFfinwin,
\qquad
F(a)=A_F(a)-I_F(a)
 \EqByDef \sum_{q\le Q_F}\left(\WinqF\right)c_q(a)-\sum_{d|a}\Irr^{(P_F)}_d F,
\quad
\forall a\in \N
$$ 
\par
\noindent
where we call $A_F$ the Analytic Part not by chance: it's, as a function of $a\in \C$, a holomorphic function, since it's a finite linear combination (not depending on $a$ \& this is important), of exponentials $e^{2\pi ija/q}$ (inside Ramanujan sums). Instead, for $F\in \AFsmwin\backslash \AFfinwin$, the regular/irregular decomposition is:\enspace $F=S_F-I_F$, 
$$
\forall F\in \AFsmwin\backslash \AFfinwin,
\qquad
F(a)=S_F(a)-I_F(a)
 \EqByDef \sum_{q\in (P_F)}\left(\WinqF\right)c_q(a)-\sum_{d|a}\Irr^{(P_F)}_d F,
\quad
\forall a\in \N
$$ 
\par
\noindent
where now $S_F\neq A_F$ is still a linear combination of exponentials above, BUT the LENGTH of $q-$sum, now, IS NOT FIXED (by Wintner's exact range, as before) : it DEPENDS ON $a$, thanks to Ramanujan Vertical Limit (in particular, depending on $p-$adic valuations of our $a$). This is new, w.r.t. [C4]. 
\smallskip
\par
\noindent
All of our Corollaries, 2.1 to 2.4, are, directly or indirectly, consequences of Theorem 2.1, namely, Wintner's Orthogonal Decomposition for the Eratosthenes Transform. All are about the, say, \lq \lq global coefficients\rq \rq, both $\WinqF$ and $\CarqF$.
 
\medskip

\par
We come, now, to $P-$local expansions, not for all arithmetic functions, like we did in Theorem 1.1 above, but in particular subsets; because, introducing mild assumptions for $F$ to expand, we want to derive a good piece of informations, for example about $P-$local coefficients. In particular, the $P-$expansion's uniqueness. 
\smallskip
\par
\noindent
Here, we give only the necessary definitions we need for next statements, leaving to $\S3$ other details. 
\par
We recall the {\bf square-free kernel}: \enspace $\kappa(1)\defineq 1$, \enspace $\kappa(n)\defineq \prod_{p|n}p$, \enspace $\forall n\in \N$; in $n=1$ these are consistent, by the convention: {\bf empty products are} $\defineq 1$. We need this kernel in order to define any arithmetic function $F$ that \lq \lq {\stampatello Ignores Prime Powers}\rq \rq, abbreviated $\IPP$, 
$$
F \enspace \IPP
\quad \definiz \quad
\forall a\in\N, \enspace F(a)\defineq F(\kappa(a)). 
$$
\par
\noindent
The standard example is: $\omega$ $\IPP$ (the prime-divisors function we recalled above). A moment's reflection brings : \enspace $F$ $\IPP$ \thinspace $\Leftrightarrow$ \thinspace $F'=\mu^2 F'$ (i.e., the Eratosthenes Transform is square-free supported), compare [C5]. 
\vfill
\par
Recall the {\bf identity function}: \enspace $\Ide(n)\defineq n$, $\forall n\in \N$; and the {\stampatello Ramanujan Conjecture}: $F(n)\ll_{\varepsilon} n^{\varepsilon}$, as $n\to \infty$, compare [C5], also for standard $\ll-$notation; while, $F$ is \lq \lq {\stampatello Neatly Sub-Linear}\rq \rq, abbreviated 
$$
F\enspace {\it is} \enspace \NSL
\quad \definiz \quad
\exists \delta<1, F(a)\ll_{\delta} a^\delta, \enspace a\to \infty. 
$$
\par
\noindent
While the $\IPP$ subset has a more arithmetic flavor, so-to-speak, the $\NSL$ is more towards analytic worlds. Next results will study both environments and their intersection, too. 
 
\eject

\par				
\noindent
In order to give next results, we define, once fixed $P\in \Primes$, {\stampatello the $q-$th $P-$Carmichael coefficient} of $F$, 
$$
\CarPqF \defineq {1\over {\varphi(q)}}\lim_x {1\over x}\sum_{a\le x}\left(\sum_{{d\in (P)}\atop {d|a}}F'(d)\right)c_q(a),
\qquad
\forall q\in \N, 
$$
\par
\noindent
of course whenever, fixed $q\in \N$, this limit exists in $\C$. Defining {\stampatello the $(P)-$restriction} of our $F$ as [C2] 
$$
F_{(P)}(a)\defineq \sum_{{d\in (P)}\atop {d|a}}F'(d),
\qquad
\forall a\in \N, 
$$
\par
\noindent
we get at once: \enspace $\CarPqF \EqByDef \Car_q\,F_{(P)}$, \enspace $\forall P\in \Primes$, $\forall q\in \N$. Next, {\stampatello the $P-$Carmichael Transform} of $F$,
$$
\CarTP F : q\in \N \mapsto \CarPqF\in \C. 
$$
\par
\noindent
Likewise, we define {\stampatello the $q-$th $P-$Wintner coefficient} of $F:\N \rightarrow \C$, if following series converges in $\C$ 
$$
\WinPqF \defineq \sum_{{d\in (P)}\atop {d\equiv 0(\!\!\bmod q)}}{{F'(d)}\over d},
\qquad
\forall q\in \N.
$$
\par
\noindent
As above,\thinspace  $\WinPqF \EqByDef \Win_q\,F_{(P)}$, \thinspace $\forall P\in \Primes$, $\forall q\in \N$. We define, when all these coefficients exist, {\stampatello the $P-$Wintner Transform} of $F$, $\WinTP F : q\in \N \mapsto \WinPqF\in \C$.  More details in $\S3$. 
\par
\noindent
\par
We state here, and prove in next subsection, following result.
\smallskip
\par
\noindent{\bf Theorem 2.2}.
%
%
 {\stampatello (Formula for $\NSL-$functions' P-Carmichael/P-Wintner coefficients)} 
\par
\noindent
{\it Let } $F$ {\it be } $\NSL$. {\it Then, } $\forall P\in \Primes$,  
\smallskip
\item{(1)} $\exists \WinPqF\in \C$, $\forall q\in \N$ \enspace {\stampatello and} \enspace  $\WinPqF=0$, $\forall q\not \in (P)$; 
\smallskip
\item{(2)} $\exists \CarPqF\in \C$, $\forall q\in \N$ \enspace {\stampatello and} \enspace  $\CarPqF=\WinPqF$, $\forall q\in \N$; 
\smallskip
\item{(3)} $\CarPqF = {\displaystyle {1\over {\varphi(q)\prod_{p\le P}\left(1-{1\over p}\right)^{-1}}}\sum_{t\in (P)}{{F(t)}\over t}c_q(t) }$, $\forall q\in (P)$; 
\smallskip
\item{(4)} $F_{(P)}(a)=F(a_{(P)})={\displaystyle \sum_{q\in (P)}(\WinPqF)c_q(a) }$, $\forall a\in \N$; 
\smallskip
\item{(5)} {\stampatello If } $G_P:\N\rightarrow \C$ {\stampatello with } $\supporto(G_P)\subseteq (P)$ {\stampatello has, for some } $\varepsilon>0$, $G_P(q)\ll_{\varepsilon,P} q^{-\varepsilon}$, $\forall q\in (P)$, {\stampatello then} 
$$
F_{(P)}(a)=\sum_{\ell\in (P)}G_P(\ell)c_{\ell}(a),\enspace \forall a\in \N
\quad \Longrightarrow \quad
G_P = \WinTP F.
$$
\par
\noindent
In other words, $\NSL$ arithmetic functions $F$ have ALL the, say, {\stampatello $(P)-$smooth Carmichael-Wintner coefficients} which expand $P-$locally our $F$, being, also, the unique such coefficients with a good decay. 
\par
We state here, and prove in next subsection, following result, too.
\smallskip
\par
\noindent{\bf Theorem 2.3}.
%
%
 {\stampatello (Formula for $\IPP-$functions' P-Carmichael/P-Wintner coefficients)} 
\par
\noindent
{\it Let } $F$ $\IPP$. {\it Then, } $\forall P\in \Primes$,  
\smallskip
\item{(1)} $\exists \WinPqF\in \C$, $\forall q\in \N$ \enspace {\stampatello and} \enspace  $\WinPqF=0$, $\forall q\not \in \Pb$; 
\smallskip
\item{(2)} $\exists \CarPqF\in \C$, $\forall q\in \N$ \enspace {\stampatello and} \enspace  $\CarPqF=\WinPqF$, $\forall q\in \N$; 
\smallskip
\item{(3)} $\CarPqF = {\displaystyle {1\over {\varphi(q)\prod_{p\le P}\left(1-{1\over p}\right)^{-1}}}\sum_{t\in (P)}{{F(t)}\over t}\cdot \left( {{\mu^2(t)t}\over {\varphi(t)}}\right)\cdot c_q(t) }$, $\forall q\in \Pb$; 
\smallskip
\item{(4)} $F_{(P)}(a)=F(a_{(P)})={\displaystyle \sum_{q\in (P)}(\WinPqF)c_q(a) }$, $\forall a\in \N$; 
\smallskip
\item{(5)} {\stampatello If } $G_P:\N\rightarrow \C$ {\stampatello has } $\supporto(G_P)\subseteq \Pb$, {\stampatello then} 
$$
F_{(P)}(a)=\sum_{\ell\in (P)}G_P(\ell)c_{\ell}(a),\enspace \forall a\in \N
\quad \Longrightarrow \quad
G_P = \WinTP F.
$$
\par
\noindent
In other words, $F$ $\IPP$ have ALL the, say, {\stampatello $(P)-$smooth Carmichael-Wintner coefficients} which expand $P-$locally our $F$, and, also, the unique such coefficients being square-free supported. 

\vfill
\eject

\par				
\noindent
{\stampatello Theorem 2.2 \& Theorem 2.3 entail the}
\smallskip
\par
\noindent{\bf Corollary 2.5}.
%
%
 {\stampatello (P-Local expansions \& their Coefficients for functions $\IPP$ and $\NSL$ too)} 
\par
\noindent
{\it Let } $F$ $\IPP$ {\it and let it be } $\NSL$, {\it too}. {\it Then, } $\forall P\in \Primes$,  
\smallskip
\item{$(i)$} $\exists \WinTP\!\!F$, $\supporto(\WinTP\!\!F)\subseteq \Pb$, $\exists \CarTP\!\!F$ {\stampatello and} \enspace  $\CarTP\!\!F=\WinTP\!\!F$; 
\smallskip
\item{$(ii)$} $\exists \WinTP\!\!\!\left(F\!\cdot\!{\Ide\over {\varphi}}\!\cdot\!\mu^2\right)$, $\exists \CarTP\!\!\!\left(F\!\cdot\!{\Ide\over {\varphi}}\!\cdot\!\mu^2\right)$ {\stampatello and} \enspace  $\CarTP\!\!\!\left(F\!\cdot\!{\Ide\over {\varphi}}\!\cdot\!\mu^2\right)=\WinTP\!\!\!\left(F\!\cdot\!{\Ide\over {\varphi}}\!\cdot\!\mu^2\right)$; 
\smallskip
\item{$(iii)$} $\CarPqF = \Car^{(P)}_q\!\!\left(F\!\cdot\!{\Ide\over {\varphi}}\!\cdot\!\mu^2\right)$, $\forall q\in \Pb$; 
\smallskip
\item{$(iv)$} $F_{(P)}(a)=F(a_{(P)})={\displaystyle \sum_{q\in (P)}(\WinPqF)c_q(a) }$, $\forall a\in \N$; 
\smallskip
\item{$(v)$} {\stampatello If } $G_P:\N\rightarrow \C$ {\stampatello with} $\supporto(G_P)\subseteq (P)$ {\stampatello has at least one, or both, of the following two properties}: 
$$
G_P = \mu^2\cdot G_P
\quad
\hbox{\stampatello or}
\quad
\exists \varepsilon>0, G_P(q)\ll_{\varepsilon,P} q^{-\varepsilon},\enspace \forall q\in (P), 
$$
\par
{\stampatello then} 
$$
F_{(P)}(a)=\sum_{\ell\in (P)}G_P(\ell)c_{\ell}(a),\enspace \forall a\in \N
\quad \Longrightarrow \quad
G_P = \WinTP F.
$$
\par
\noindent{\bf Proof}. Since $(iii)$ is the \lq \lq {\stampatello marriage}\rq \rq, so to speak, {\stampatello of Th.m 2.2's $(3)$ and Th.m 2.3's $(3)$}, we only prove $(ii)$. 
\par
\noindent
For this, we only need to show: 
$$
F\; \NSL 
\enspace \Rightarrow \enspace 
F\cdot {\Ide\over {\varphi}}\cdot \mu^2\; \NSL 
\leqno{(\ast)}
$$
\par
\noindent
because, from this, we may apply {\stampatello Th.m 2.2's first part of $(1)$ and $(2)$ for \enspace $F\!\cdot\!{\Ide\over {\varphi}}\!\cdot\!\mu^2$\enspace getting} $(ii)$, in this way. In turn, our $(\ast)$ above follows from, when $n\to \infty$, 
$$
\mu^2(n)=1
\enspace \Rightarrow \enspace
{{\Ide(n)}\over {\varphi(n)}}={n\over {\varphi(n)}}
 =\prod_{p|n}{p\over {p-1}}
  =\prod_{p|n}\left(1+{1\over {p-1}}\right)
   =\exp\left(\sum_{p|n}\log\left(1+{1\over {p-1}}\right)\right)
    \ll \log n. 
\leqno{(\ast\ast)}
$$
\par
\noindent
This bound, in fact, comes from 
$$
\exp\left(\sum_{p|n}\log\left(1+{1\over {p-1}}\right)\right)\ll \exp\left(\sum_{p|n}{1\over {p-1}}\right)
 = \exp\left(\sum_{p|n}{1\over p}+\sum_{p|n}{1\over {p(p-1)}}\right)
  \ll \exp\left(\sum_{p|n}{1\over p}\right); 
$$
\par
\noindent
the trivial bound \enspace $\sum_{p|n}{1\over p}\ll \sum_{p\le n}{1\over p}$ \enspace with (see [T])
$$
\sum_{p\le x}{1\over p}\ll \log \log x\qquad
\hbox{\rm as}\enspace x\to \infty,
$$
\par
\noindent
thus give $(\ast\ast)$.\hfill $\square$

\bigskip

\par
\noindent
In the following, we'll write QED for the end of a Proof's part. The Proof, except for Facts, ends with \hfill $\square$

\medskip

\par
\noindent
Hereafter we define $\widetilde{F}$, {\stampatello the \lq \lq IPPification\rq \rq}, say, of an arithmetic function $F$ as
$$
\forall a\in \N,
\quad
\widetilde{F}(a)\defineq F(\kappa(a)). 
$$
\par
\noindent
Whatever is $F:\N \rightarrow \C$, its IPPification \enspace $\widetilde{F}$ \enspace is $\IPP$ by definition; also, $F$ $\IPP$ entails $\widetilde{F}=F$.  

\vfill
\eject

\par				
\noindent
{\stampatello Theorem 2.2 \& Theorem 2.3 also entail the other}
\smallskip
\par
\noindent{\bf Corollary 2.6}.
%
%
 {\stampatello (P-Local expansions for $\NSL$ functions $F$ and $\mu^2 F\Ide/\varphi$)} 
\par
\noindent
{\it Let } $F$ {\it be } $\NSL$. {\it Fix } $P\in \Primes$. {\it Then}  
\smallskip
\item{$(i)$} {\it Carmichael $P-$Transforms and Wintner $P-$Transforms} {\stampatello exist and are equal, for the $3$ functions } $F$, $\widetilde{F}$ {\stampatello and} $\mu^2\cdot F\cdot{{\Ide}\over {\varphi}}$; 
\smallskip
\item{$(ii)$} $q\in \Pb$ \quad $\Rightarrow$ \quad $\Car^{(P)}_q\widetilde{F}=\Car^{(P)}_q\!\left(F\!\cdot\!{\Ide\over {\varphi}}\!\cdot\!\mu^2\right)$; 
\smallskip
\item{$(iii)$} ${\displaystyle \left({a\over {\varphi(a)}}-1\right) F(a) = \sum_{{\ell\in(P)}\atop {{\ell \,  \hbox{\piccolissimo CUBE-FREE,\thinspace NOT\thinspace \thinspace S-F}}\atop {\left.{{\ell}\over {\kappa(\ell)}}\right|a}}}\Win^{(P)}_{\ell}\!\left(\mu^2 F\cdot {{\Ide}\over {\varphi}}\right)c_{\ell}(a), \quad \forall a\in \Pb }$.
\medskip
\par
\noindent{\bf Remark 1}.
%
%
 We use now: when $F$ is $\NSL$, then $G$ has {\stampatello Ramanujan Conjecture} \enspace $\Rightarrow$ \enspace $FG$ is $\NSL$.\hfill $\diamond$ 
\medskip
\par
\noindent{\bf Proof}. Using $(1)$ and $(2)$ of {\stampatello Theorem} 2.2, since: $F$ $\NSL$ $\Rightarrow$ $\widetilde{F}$ $\NSL$ is trivial, to prove $(i)$ we only need the property $(\ast)$ in previous Corollary 2.5 Proof.\hfill QED 
\par
\noindent
Since $\widetilde{F}$ $\IPP$ by definition, present $(ii)$ comes from $(iii)$ of Corollary 2.5.\hfill QED
\par
\noindent
The {\stampatello interesting part is formula $(iii)$, for which we start from two formul\ae:}
\par
\noindent
{\stampatello using $(1)$, $(2)$, $(4)$ of Theorem 2.2} for $\mu^2 F\Ide /\varphi$ and $(\ast)$ above, we get 
$$
F \enspace \NSL \enspace \hbox{\stampatello and} \enspace a\in \Pb
\quad \Longrightarrow \quad
F(a){a\over {\varphi(a)}}=\sum_{\ell\in (P)}\Car^{(P)}_{\ell}\!\left(\mu^2 F\!\cdot\!{\Ide\over {\varphi}}\right)c_{\ell}(a)
\leqno{(I)}
$$
\par
\noindent
{\stampatello and by Theorem 2.3's $(1)$, $(2)$, $(3)$, $(4)$ \& Theorem 2.2's $(3)$}, we obtain, recalling $\sumflat_{\ell \cdots}\EqByDef \sum_{\ell \cdots, \ell \, \hbox{\piccolo S-F}}$ , 
$$
\widetilde{F} \thinspace \IPP \enspace \hbox{\stampatello and} \enspace \NSL, \enspace a\in \Pb
\quad \Longrightarrow \quad
\widetilde{F}(a)=\sumflat_{\ell\in (P)}\Car^{(P)}_{\ell}\!\left(\mu^2 F\!\cdot\!{\Ide\over {\varphi}}\right)c_{\ell}(a)
\leqno{(II)}
$$
\par
\noindent
whence, using \enspace $\widetilde{F}\cdot \mu^2=F\cdot \mu^2$, 
$$
a\in \Pb
\quad \Longrightarrow \quad
F(a)\buildrel{(II)}\over{=\!=\!=}\sum_{\ell\in (P)}\Car^{(P)}_{\ell}\!\left(\mu^2 F\!\cdot\!{\Ide\over {\varphi}}\right)c_{\ell}(a)-\sum_{{\ell\in (P)}\atop {\ell\, \hbox{\piccolissimo NOT\thinspace S-F}}}\Car^{(P)}_{\ell}\!\left(\mu^2 F\!\cdot\!{\Ide\over {\varphi}}\right)c_{\ell}(a)
\buildrel{(I)}\over{=\!=\!=}
$$
$$
\buildrel{(I)}\over{=\!=\!=} F(a){a\over {\varphi(a)}}-\sum_{{\ell\in (P)}\atop {\ell\, \hbox{\piccolissimo NOT\thinspace S-F}}}\Win^{(P)}_{\ell}\!\left(\mu^2 F\!\cdot\!{\Ide\over {\varphi}}\right)c_{\ell}(a), 
$$
\par
\noindent
recalling {\stampatello Th.m} 2.2's $(2)$ for $\mu^2 F\!\cdot\!{\Ide\over {\varphi}}$. 
\par
\noindent
{\stampatello Thus}
$$
\left({a\over {\varphi(a)}}-1\right)F(a)=\sum_{{\ell\in (P)}\atop {\ell\, \hbox{\piccolissimo NOT\thinspace S-F}}}\Win^{(P)}_{\ell}\!\left(\mu^2 F\!\cdot\!{\Ide\over {\varphi}}\right)c_{\ell}(a),
\quad \forall a\in \Pb. 
$$
\par
\noindent
The {\stampatello supplementary conditions}: $\ell$ \enspace {\stampatello cube-free and} \enspace ${\displaystyle \left.{{\ell}\over {\kappa(\ell)}}\right|a }$ \enspace are needed, just {\stampatello to avoid} $c_{\ell}(a)=0$: 
$$
c_{\ell}(a)\buildrel{\hbox{\rm [K]}}\over{=\!=\!=}\sum_{{d|a}\atop {d|\ell}}d\mu\left({{\ell}\over d}\right)
 \buildrel{\mu^2(a)=1}\over{=\!=\!=\!=\!=}\sum_{{d|a}\atop {d|\kappa(\ell)}}d\mu\left({{\ell}\over {\kappa(\ell)}}\cdot {{\kappa(\ell)}\over d}\right)
  = \mu\left({{\ell}\over {\kappa(\ell)}}\right)\sum_{{d|a,d|\kappa(\ell)}\atop {\left({{\ell}\over {\kappa(\ell)}},{{\kappa(\ell)}\over d}\right)=1}}d\mu\left({{\kappa(\ell)}\over d}\right)\neq 0
$$
\par
\noindent
{\stampatello entails}, in particular, $\mu(\ell/\kappa(\ell))\neq 0$ $\Rightarrow$ $\ell/\kappa(\ell)$ {\stampatello square-free} $\Rightarrow$ $\ell$ is {\stampatello cube-free}\hfill (compare the following) 
\par
\noindent
{\stampatello and the inner $d-$sum has}:
$$
\left({{\ell}\over {\kappa(\ell)}},{{\kappa(\ell)}\over d}\right)=1
\enspace \Longleftrightarrow \enspace 
\left( \left.p\right|{{\ell}\over {\kappa(\ell)}} \thinspace \Rightarrow \thinspace p|d \right)
\enspace \Longleftrightarrow \enspace 
\left.{{\ell}\over {\kappa(\ell)}}\right|d, 
$$
\par
\noindent
using in this last equivalence (compare implications above): $\ell$ cube-free $\Leftrightarrow$ $\ell/\kappa(\ell)$ square-free; finally, last condition above, from $d|a$, entails ${{\ell}\over {\kappa(\ell)}}|a$.\hfill $\square$
\medskip
\par
See that also Lemma 4.2 and Remark 5 entail, in particular, the cube-free restriction in the $\ell-$sum of $(iii)$, above. 
\vfill
\eject

\par				
\noindent{\bf 2.1 Theorems' Proofs requiring sections 3,4 Lemmata} 
\bigskip
\par
\noindent
Hereafter we use the $O-$notation, a perfect synonym of $\ll$, including subscripts: see [D],[T] \& compare [C5]. 
\medskip
\par
\noindent
Our Lemmata $3.1$, $3.4$, $3.6$, $3.7$ and $4.1$ prove {\bf Theorem 2.2}. 
\par
\noindent{\bf Proof}. Second part of $(1)$ is trivial from $\WinPqF$ definition, while first part follows from $F\; \NSL$ and formula $(1P)$ of Lemma $3.6$ giving even more than convergence, absolute convergence: 
$$
\WinPqF \EqByDef \sum_{{d\in (P)}\atop {d\equiv 0(\bmod q)}}{{F'(d)}\over d}\ll_{\delta} q^{\delta-1}\sum_{m\in (P)}m^{\delta-1}
 \ll_{\delta,P} q^{\delta-1}.
$$
\enspace \hfill QED 
\par
\noindent
Then, Lemma $3.7$ proves immediately $(2)$, from $F\; \NSL$ $\Rightarrow$ $F'\; \NSL$ and $(2P)$ of Lemma 3.6.\hfill QED 
\par
\noindent
For $(3)$, $\CarPqF$ definition and Lemma $3.4$ {\stampatello entail}, $\forall P\in \Primes$ \& $\forall q\in (P)$, 
$$
\CarPqF = {1\over {\varphi(q)}}\lim_x {1\over x}\sum_{a\le x}c_q(a)\sum_{{t\in (P)}\atop {{t|a}\atop {{a\over t}\in )P(}}}F(t)
         = {1\over {\varphi(q)}}\lim_x {1\over x}\sum_{{t\in (P)}\atop {t\le x}}F(t)c_q(t)\sum_{{m\le {x\over t}}\atop {m\in )P(}}1, 
$$
\par
\noindent
{\stampatello from}: $c_q(tm)=c_q(t)$, $\forall q\in (P)$, $\forall m\in )P($; now, $(3P)$ of Lemma $3.6$, with $x/t$ instead of $x$, gives: 
$$
\CarPqF = {1\over {\varphi(q)}}\lim_x {1\over x}\sum_{{t\in (P)}\atop {t\le x}}F(t)c_q(t)\left(\prod_{p\le P}\left(1-{1\over p}\right){x\over t}+O_P(1)\right)
$$
$$
= {1\over {\varphi(q)\prod_{p\le P}\left(1-{1\over p}\right)^{-1}}}\lim_x \sum_{{t\in (P)}\atop {t\le x}}{{F(t)}\over t}c_q(t),
$$
\par
\noindent
because:  
$$
{1\over {\varphi(q)}}\lim_x {1\over x}\sum_{{t\in (P)}\atop {t\le x}}F(t)c_q(t)O_P(1)
= O_{P,\delta}\left(\lim_x {1\over x}\sum_{{t\in (P)}\atop {t\le x}}t^{\delta}\right)
 = O_{P,\delta,\varepsilon}\left(\lim_x {1\over x}x^{\delta}\cdot x^{\varepsilon}\right),
$$
\par
\noindent
{\stampatello from}: the trivial $|c_q(t)|\le \varphi(q)$, $\forall t\in \N$ and $(2P)$ of Lemma $3.6$, so this {\stampatello vanishes choosing} $0<\varepsilon<1-\delta$; {\stampatello main term's series over $t\in (P)$ converges,} even better, absolutely: 
$$
\sum_{{t\in (P)}\atop {t\le x}}{{|F(t)|}\over t}|c_q(t)|\ll_{q,\delta} \sum_{t\in (P)}t^{\delta-1}
 \ll_{q,\delta,P} 1, 
$$
\par
\noindent
again from $(1P)$ quoted above.\hfill QED
\par
\noindent
The {\stampatello local expansion} in $(4)$ comes straight from $(1)$ and Lemma $3.1$.\hfill QED
\par
\noindent
To show $(5)$: 
$$
\ell\in (P)
\buildrel{(3)}\over{\Longrightarrow} 
\Car^{(P)}_{\ell}\,F = {1\over {\varphi(\ell)\prod_{p\le P}\left(1-{1\over p}\right)^{-1}}}\sum_{t\in (P)}{{c_{\ell}(t)}\over t}\sum_{q\in (P)}G_P(q)c_q(t)
 \buildrel{\hbox{\stampatello Lemma\thinspace 4.1}}\over{=\!=\!=\!=\!=\!=\!=\!=\!=} G_P(\ell), 
$$
\par
\noindent
{\stampatello provided } $\sum_t$ {\stampatello and} $\sum_q$ {\stampatello may be exchanged; this} follows {\stampatello from the $t,q$ double series absolute convergence}, as:
$$
\sum_{t\in (P)}{{|c_{\ell}(t)|}\over t}\sum_{q\in (P)}|G_P(q)|\cdot |c_q(t)|
\ll_{\ell,\varepsilon,P,q} \sum_{t\in (P)}t^{-1}\sum_{q\in (P)}q^{-\varepsilon}
 \ll_{\ell,\varepsilon,P,q} 1,
$$
\par
\noindent
{\stampatello thanks to} $(1P)$ one last time.\hfill $\square$ 

\vfill
\eject

\par				
\noindent
Recall \enspace $\Pb \EqByDef (P)\cap \{n\in \N : n$ square-free$\}$. 
\par
\noindent
Our Lemmata $3.1$, $3.3$, $3.4$, $3.6$, $3.7$ and $4.1$ altogether prove {\bf Theorem 2.3}.
\par
\noindent{\bf Proof}. We {\stampatello closely follow Th.m 2.2 Proof, for properties} $(1)$, $(2)$ \& $(4)$. 
\par
\noindent
Second part of $(1)$ is trivial from $\WinPqF$ definition; this, with $F'=\mu^2\cdot F'$, gives first part, too.\hfill QED 
\par
\noindent
Lemma $3.7$ \& $F$ $\IPP$, again from $F'=\mu^2\cdot F'$, immediately give $(2)$.\hfill QED 
\par
\noindent
For $(3)$, Th.m 2.2 Proof applies here, once {\stampatello fixed} $P\in \Primes$ and $q\in (P)$, {\stampatello up to equation}
$$
\CarPqF = {1\over {\varphi(q)}}\lim_x {1\over x}\sum_{{t\in (P)}\atop {t\le x}}F(t)c_q(t)\left(\prod_{p\le P}\left(1-{1\over p}\right){x\over t}+O_P(1)\right); 
$$
\par
\noindent
here, hypothesis $F$ $\IPP$ gives, trivially, $F=\widetilde{F}$ : we apply Lemma $3.3$ to $F(t)=\widetilde{F}(t)$, $\forall t\in \N$, {\stampatello to get}
$$
{1\over {\varphi(q)}}\lim_x {1\over x}\sum_{{t\in (P)}\atop {t\le x}}F(t)c_q(t)\prod_{p\le P}\left(1-{1\over p}\right){x\over t}
={1\over {\varphi(q)\prod_{p\le P}\left(1-{1\over p}\right)^{-1}}}\lim_x \sum_{{t\in (P)}\atop {t\le x}}{{c_q(t)}\over t}\sum_{{d|t}\atop {\kappa({t\over d})|d}}\mu^2(d)F(d) 
$$
\par
\noindent
and {\stampatello setting} $t=dm$ this double sum becomes, since $d\in (P)$ and $\kappa(m)|d$ entail $m\in (P)$, the following: 
$$
\sum_{{d\in (P)}\atop {d\le x}}{{\mu^2(d)F(d)}\over d}\sum_{{m\in (P)}\atop {{m\le x/d}\atop {\kappa(m)|d}}}{{c_q(dm)}\over m}
=\sum_{{d\in (P)}\atop {d\le x}}{{\mu^2(d)F(d)}\over d}\cdot c_q(d)\sum_{{m\le x/d}\atop {\kappa(m)|d}}{1\over m},
\quad \forall q\in \Pb, 
$$
\par
\noindent
because $\mu^2(q)=1$ $\Rightarrow$ (see $\S3.1.2$) $\Rightarrow$ $c_q(a)$ $\IPP$ as a function of $a\in \N$ $\Rightarrow$ $c_q(dm)=c_q(d)$, {\stampatello for any} $m\in \N$ {\stampatello with} $\kappa(m)|d$, because in this case $p|m$ $\Rightarrow$ $p|d$; as soon as $x\ge {\displaystyle \prod_{p\le P}p }$, we have $d\in \Pb$ $\Rightarrow$ $x\ge d$; {\stampatello thus}
$$
\lim_x \sum_{{t\in (P)}\atop {t\le x}}{{c_q(t)}\over t}\sum_{{d|t}\atop {\kappa({t\over d})|d}}\mu^2(d)F(d)
=\sum_{d\in (P)}{{\mu^2(d)F(d)}\over d}\cdot c_q(d)\cdot \lim_x \sum_{{m\le x/d}\atop {\kappa(m)|d}}{1\over m}, 
$$
\par
\noindent
where 
$$
\lim_x \sum_{{m\le x/d}\atop {\kappa(m)|d}}{1\over m}=\sum_{\kappa(m)|d}{1\over m}
 = \prod_{p|d}\sum_{K=0}^{\infty}p^{-K}
  = \prod_{p|d}{1\over {1-{1\over p}}}
   = \left(\prod_{p|d}\left(1-{1\over p}\right)\right)^{-1}
    = \left({{\varphi(d)}\over d}\right)^{-1}
$$
\par
\noindent
{\stampatello gives} $(3)$'s RHS(Right Hand Side), once we prove that {\stampatello the remainder} $O_P(1)$ {\stampatello in equation above has $\lim_x$ vanishing} : this is from
$$
{1\over {\varphi(q)}}\, {1\over x}\sum_{{t\in (P)}\atop {t\le x}}F(t)c_q(t)O_P(1)
\ll_{q,P} {1\over x}\sum_{{t\in (P)}\atop {t\le x}}\sum_{{d|t}\atop {\kappa({t\over d})|d}}\mu^2(d)|F(d)| 
 \ll_{q,P} {1\over x}\sum_{d\in \Pb}|F(d)| \sum_{{m\in (P)}\atop {m\le x/d}}1
  \ll_{q,P,\varepsilon,F} x^{\varepsilon-1}
   \buildrel{x}\over{\to} 0,
$$
\par
\noindent
{\stampatello using} $(2P)$ of Lemma $3.6$.\hfill QED 
\par
\noindent
For $(4)$, use $(1)$ with Lemma $3.1$.\hfill QED 
\par
\noindent
For $(5)$, on the same line of {\stampatello Th.m} 2.2 {\stampatello Proof, this time we use $G_P = \mu^2 \cdot G_P$ to get the absolute convergence}.\hfill $\square$

\vfill
\eject

\par				
\noindent{\bf 3. Technical Lemmas} 
\bigskip
\par
\noindent
\par
\noindent{\bf 3.1 Ramanujan vertical limit: non-vanishing Ramanujan sums} 
\medskip
\par
\noindent
We start with a small Lemma, a {\stampatello Fact}, regarding Ramanujan Sums. We use it at once, in next Lemma $3.1$. 
\smallskip
\par
\noindent{\bf Fact 3.1}.
%
%
 {\stampatello (Ramanujan vertical limit)} 
$$
c_q(a)\neq 0
\enspace \Longrightarrow \enspace
v_p(q)\le v_p(a)+1,
\enspace
\forall p|q.
\leqno{\Rvl}
$$
\par
\noindent{\bf Remark 2}.
%
%
 Actually, RHS inequality above holds $\forall p\in \Primes$ (trivially: $\forall p\in \Primes$ fixed, $v_p(q)\ge 0$, $\forall q\in \N$).\hfill $\diamond$  
\medskip

\par
In fact, we prove even more than {\it Ramanujan vertical limit} above, since we have the following equivalence. 
\smallskip
\par
\noindent{\bf Proposition 3.1}.
%
%
 {\stampatello (Non-vanishing Ramanujan sums)} 
$$
c_q(a)\neq 0
\enspace \Longleftrightarrow \enspace
v_p(q)\le v_p(a)+1,
\enspace
\forall p\in \Primes.
$$
\par
\noindent{\bf Proof}. We apply H\"older's formula: (compare $(1.6)$, Fact 1.5 in [C5]) 
$$
c_q(a)=\varphi(q)\mu(q/(q,a))/\varphi(q/(q,a)),
\quad \forall q\in \N, \forall a \in \Z, 
$$
\par
\noindent
to conclude that the Ramanujan sum $c_q(a)$ doesn't vanish IFF $\mu(q/(q,a))\neq 0$; this, in turn, is equivalent to: $q/(q,a)$ is square-free, so, 
$$
c_q(a)\neq 0
\thinspace \Leftrightarrow \thinspace
v_p(q/(q,a))\le 1, \forall p\in \Primes
\thinspace \Leftrightarrow \thinspace
v_p(q)-v_p((q,a))\le 1, \forall p\in \Primes, 
$$
\par
\noindent
from the complete additivity of $v_p$, $\forall p\in \Primes$; then, by definition $v_p((q,a))=\min(v_p(q),v_p(a))$, whence 
$$
c_q(a)\neq 0
\thinspace \Leftrightarrow \thinspace
v_p(q)-\min(v_p(q),v_p(a))\le 1, \forall p\in \Primes,
\thinspace \Leftrightarrow \thinspace
\max(0,v_p(q)-v_p(a))\le 1, \forall p\in \Primes;
$$
\par
\noindent
we've two possible cases for the maximum, for each fixed $p\in \Primes$ : first case \enspace $v_p(q)\le v_p(a)$, with vanishing $\max$, renders the RHS (right hand side) the trivial $0\le 1$; and second case \enspace $v_p(q)>v_p(a)$, having $\max=v_p(q)$, entails that RHS is equivalent to \enspace $v_p(q)-v_p(a)\le 1$, i.e. $v_p(q)\le v_p(a)+1$.\hfill $\square$  

\bigskip

\par
\noindent{\bf 3.1.1 Local expansions: horizontal limits on moduli} 
\medskip
\par
\noindent
Once fixed $P\in \Primes$, the existence of {\stampatello Wintner $P-$smooth Transform} entails following formula. 
\smallskip
\par
\noindent{\bf Lemma 3.1}.
%
%
 {\stampatello (Wintner $P-$smooth Transform existence gives $P-$local expansion)}
\par
\noindent
{\it Let } $F:\N \rightarrow \C$ {\it have, } {\stampatello for a fixed } $P\in \Primes$, {\it the Wintner $P-$smooth Transform. Then}
$$
F(a)=\sum_{q\in (P)}\left( \WinPqF \right)c_q(a), 
\quad
\forall a\in (P). 
$$
\smallskip
\par
\noindent{\bf Proof}. Fix $P\in \Primes$ and choose $a\in (P)$ : $\Rvl$ $\Rightarrow $ next ${\displaystyle \sum_{q\in (P)}}$ is finite, whence $x-$independent, as $x\to \infty$, 
$$
F(a)=\sum_{{d\in (P)}\atop {d|a}}F'(d)
 \buildrel{(1.1)}\over{=\!=\!=}\lim_x \sum_{{d\in (P)}\atop {d\le x}}{{F'(d)}\over d}\sum_{q|d}c_q(a)
  \buildrel{\hbox{\stampatello (Rvl)}}\over{=\!=\!=\!=\!=}\sum_{q\in (P)}\left( \WinPqF \right)c_q(a), 
$$
\par
\noindent
from $\WinPqF$ definition, $\forall q\in \N$, thanks to our existence hypothesis and the property that the set $(P)$ is divisor-closed, here: $d\in (P)$, $q|d$ $\Rightarrow $ $q\in (P)$. This $q-$sum is bounded ONLY in terms of $P$ and $a$.\hfill $\square$ 
\par
\noindent{\bf Remark 3}.
%
%
 The power of $\Rvl$ is in the exchange of limits: in above Proof, these are over $q\in (P)$ and over $x\to \infty$; actually, first one is a fake limit process, as $q$ is bounded {\stampatello uniformly } $\forall x\in \N$, by $\Rvl$ ! And we need ONLY the existence (guaranteed by existence hypothesis) of limit over $x\to \infty$.\hfill $\diamond$  

\bigskip

\par				
\noindent{\bf 3.1.2 Ramanujan sums ignoring prime powers: square-free moduli} 
\medskip
\par
\noindent
We invoke Kluyver's Formula [K] quoted above, in order to prove that: once fixed $q\in \N$ square-free, the arithmetic function
$$
c_q(\cdot): a\in \N \mapsto c_q(a)\in \C
$$
\par
\noindent
is an $\IPP$ function. In fact, $q$ square-free entails all $d|q$, its divisors, are square-free, too:
$$
c_q(a)=\sum_{{d|a}\atop {d|q}}d\mu\left({q\over d}\right)
 =\sumflat_{{d|a}\atop {d|q}}d\mu\left({q\over d}\right)
  =\sum_{{d|\kappa(a)}\atop {d|q}}d\mu\left({q\over d}\right)
   =c_q(\kappa(a)), 
$$
\par
\noindent
giving : {\bf Fact 3.2}.
%
%
 $\mu^2(q)=1$ $\Rightarrow$ $c_q(\cdot)$ $\IPP$.\hfill QED
\medskip
\par
\noindent
(Yes, we use QED also when proving Facts: recall, they're small Lemmas.)

\bigskip

\par
\noindent{\bf 3.1.3 Horizontal and vertical constraints combined: Finite-Length expansions} 
\medskip
\par
\noindent
As a corollary to above Lemma $3.1$ we get next Lemma $3.2$, using \enspace $\left|\Pb\right|=2^{\pi(P)}<\infty$, $\forall P\in \Primes$ {\stampatello fixed}; however, we give a shorter independent Proof, following. 
\smallskip
\par
\noindent{\bf Lemma 3.2}.
%
%
 {\stampatello ($P-$smooth square-free arguments always have $P-$local expansion)}
\par
\noindent
{\it Let } $F:\N \rightarrow \C$ {\stampatello and choose any } $P\in \Primes$. {\it Then, } {\stampatello recalling } 
$
{\displaystyle \WinPbqF \EqByDef \sum_{{d\in \Pb}\atop {d\equiv 0(\!\! \bmod q)}}{{F'(d)}\over d}
}, 
$
\enspace {\stampatello we have }
$$
F(a)=\sum_{q\in \Pb}\left(\WinPbqF\right)c_q(a),
\quad 
\forall a\in \Pb, 
$$
$$
F_{\Pb}(a)\defineq \sum_{{d\in \Pb}\atop {d|a}}F'(d)=\sum_{q\in \Pb}\left(\WinPbqF\right)c_q(a),
\quad 
\forall a\in \N 
$$
\par
\noindent
{\stampatello and, joining the hypothesis $F$ $\IPP$,}
$$
F_{(P)}(a)=\sum_{q\in (P)}\left(\WinPqF\right)c_q(a),
\quad 
\forall a\in \N, 
\enspace 
\hbox{\stampatello too}.
$$
\smallskip
\par
\noindent{\bf Proof}. Just one line: 
$$
a\in \Pb 
\enspace \Longrightarrow \enspace 
F(a)=\sum_{{d\in \Pb}\atop {d|a}}F'(d)
 \buildrel{(1.1)}\over{=\!=\!=}\sum_{d\in \Pb}{{F'(d)}\over d}\sum_{q|d}c_q(a)
  \buildrel{|\Pb|<\infty}\over{=\!=\!=\!=\!=}\sum_{q\in \Pb}\left( \WinPbqF \right)c_q(a). 
\enspace 
\hbox{\QED}
$$
\par
\noindent
First formula is proved, we leave other two formul\ae \enspace proofs to the interested readers.\enspace \hfill $\square$ 

\bigskip

\par
\noindent{\bf 3.1.4 Odds \& ends for horizontal and vertical limits} 
\medskip
\par
\noindent
We saw above already the horizontal limits: we are limiting the prime-factors of a natural, up to a certain \lq \lq limit\rq \rq: the classic example are the $P-$smooth naturals, whose set $(P)$ we are using a lot, in this paper. More in general we may, so-to-speak, agree that we may have any kind of horizontal limit, simply FIXING a subset of primes, for our natural's prime-factors.
\par
\noindent
We saw vertical limits, too: for example, the limitations on $p-$adic valuations in Fact 3.1 above. More in general we may, so-to-speak, complicate our life with arbitrary vertical limits, not FIXING a single number, say, $V$, for our $p-$adic valuations to be $v_p(n)\le V$, for our natural $n$. For example, asking $n$ square-free is a FIXED vertical limit: ALL $v_p(n)\le 1$, for our $n$. 
\par
\noindent
Our Theorem 1.1 and, more in detail, its Lemma 3.2 counterpart, use the fact that a combination of a vertical (fixed) limit with a horizontal one, of course, entails a finite number of naturals with both constraints: \lq \lq $p|n$ implies $p\le P$\rq \rq \enspace and \lq \lq $v_p(n)\le 1$ for all $p$\rq \rq, actually, give exactly $2^{\pi(P)}$ numbers $n$. 
\par				
\noindent
Generalizing even more, we may start from the, say, Erd\"os FACTORIZATION: $\N = )P( \otimes (P)$, meaning $a=rs$ we saw at page 3, going to more \& more general ways, with (like for previous) two factors or even with a finite ($>2$) number of \lq \lq factors\rq \rq, choosing both the primes subsets (for the \lq \lq horizontal choice\rq \rq) and the $p-$adic valuations subsets (for the \lq \lq vertical choice\rq \rq).
\par
One common feature, however, should be the ORTHOGONALITY of these sets that are the $\N-$factors. In our, say, $\WOD-$application, orthogonality is simply the COPRIMALITY, in between factors' elements.
\par
In the following Lemmata, we are actually still surfing horizontal/vertical limits sea: esp., next Lemma 3.3 is about the \lq \lq IPPification\rq \rq, that's a kind of giving to $F$'s divisors (i.e., to $F'$ support) the FIXED vertical limit 1, because (see page 8) $F$ $\IPP$ is equivalent to $F'$ square-free supported. 
\par
In the subsequent Lemma 3.4 we consider instead the $(P)-$restriction of our $F$ : a kind of restricting $F$'s divisors (i.e., $F'$ support) to the FIXED horizontal limit $P$. 
\medskip
\par
Next result is a kind of explicit formula, for the IPPification of an arbitrary $F$; even if it may appear as trivial, it will be very useful, when applied inside other formul\ae; its name, then, is underlining the universality, as it holds for all arithmetic functions. \enspace Recall abbreviations: o.w. $=$ otherwise, $\assurdo$ $=$ absurd.  
\smallskip
\par
\noindent{\bf Lemma 3.3}.
%
%
 {\stampatello (IPPification Universal Formula)}
\par
\noindent
{\it Let } $F:\N \rightarrow \C$ {\stampatello and recall } $\widetilde{F}(a)\EqByDef F((\kappa(a))$, $\forall a\in \N$. {\it Then,}
$$
\forall a\in \N, 
\quad
\widetilde{F}(a)=\sum_{t|a,\kappa({a\over t})|t}\mu^2(t)F(t).
$$
\smallskip
\par
\noindent{\bf Proof}. Just one line: 
$$
t|\kappa(a), \kappa({a\over t})\left.\right|t
\Rightarrow 
\kappa(a)=tm, \mu^2(tm)=1, (t,m)=1
\Rightarrow 
m=1:o.w.\exists p|m 
\Rightarrow 
p|a\, \&\, p\nondivide t
\Rightarrow 
p\left|\kappa({a\over t})\right.
\Rightarrow 
p|t \assurdo
$$
\enspace \hfill \square

\medskip

\par
Next result is Lemma $3.1$ in [C2]: once fixed $P\in \Primes$, for any $a\in \N$ we have the representation with a {\stampatello unique} $(r,t)\in )P( \times (P)$ \thinspace as \thinspace $a=rt$. Notice that next Lemma, too, is Universal (compare the above). 
\smallskip
\par
\noindent{\bf Lemma 3.4}.
%
%
 {\stampatello (M\"{o}bius Switch)}
\par
\noindent
{\it Let } $F:\N \rightarrow \C$ {\stampatello and choose } $P\in \Primes$. {\it Then, recalling } $F_{(P)}(a)\EqByDef F(a_{(P)})$, $\forall a\in \N$, 
$$
\forall a\in \N,
\quad 
F_{(P)}(a)=\sum_{{t\in (P)}\atop {{t|a}\atop {{a\over t}\in )P(}}}F(t).
$$
\smallskip
\par
\noindent{\bf Proof}. Another one line (recall notations $a_{(P)}$ and $a_{)P(}$ at page 3): 
$$
t\in (P), t|a, {a\over t}:=r\in )P( 
\thinspace \thinspace \Rightarrow \thinspace 
t|a_{(P)}, a=a_{(P)}\,a_{)P(}
\thinspace \Rightarrow \thinspace 
{{a_{(P)}}\over t}\,a_{)P(}=r\in )P(
\thinspace \thinspace \Rightarrow \thinspace 
{{a_{(P)}}\over t}\in (P)\cap )P(
\thinspace \thinspace \Rightarrow \thinspace 
t=a_{(P)}.
\enspace 
\square
$$
\medskip
\par
We marry these 2 Lemmas in next. With $a_{\Pb}\defineq {\displaystyle \prod_{{p\le P}\atop {p|a}}p }$, $F_{\Pb}(a)\EqByDef F(a_{\Pb})$, $\forall a\in \N$, see Lemma 3.2. 
\smallskip
\par
\noindent{\bf Lemma 3.5}.
%
%
 {\stampatello (Local $P-$expansions and $P-$smooth square-free restrictions)}
\par
\noindent
{\it Let } $F:\N \rightarrow \C$ {\stampatello and choose } $P\in \Primes$. {\it Then, recalling } $\widetilde{F}_{(P)}(a)=\widetilde{F_{(P)}}(a)=F(a_{\Pb})$, $\forall a\in \N$, 
$$
\forall a\in \N,
\quad 
F(a_{\Pb})=\sum_{
{t|a_{(P)}}\atop {\kappa({a_{(P)}\over t})|t}
}\mu^2(t)F(t)
=\sum_{
{t\in (P)}\atop 
{
{t|\kappa(a)}\atop {{{\kappa(a)}\over t}\in )P(}
}
}F(t)
=\sum_{
{t\in \Pb}\atop {
{t|a}\atop {
{{\kappa(a)}\over t}\in )P(}
}
}
F(t).
$$
\smallskip
\par
\noindent{\bf Proof}. We prove: Rightmost RHS $=$ Leftmost LHS(Left Hand Side), in two lines, as $)P($ is divisor-closed, 
$$
t\in \Pb, t|a, {{\kappa(a)}\over t}\in )P(
\thinspace \Rightarrow \thinspace 
t|a_{\Pb}
\thinspace \Rightarrow \thinspace 
t|\kappa(a_{(P)})
\thinspace \Rightarrow \thinspace 
$$
$$
\thinspace \Rightarrow \thinspace 
{{\kappa(a_{(P)})}\over t}\kappa(a_{)P(})={{\kappa(a)}\over t}\in )P(
\thinspace \Rightarrow \thinspace 
{{\kappa(a_{(P)})}\over t}\in (P)\cap )P(
\thinspace \Rightarrow \thinspace 
t=\kappa(a_{(P)})=a_{\Pb}. 
$$
\enspace \hfill \square
\par
Lemma 3.5 speaks about arguments $a\in \Pb$, but not like Lemma 3.2 that gives $P-$local expansions: it combines Lemmas 3.3 and 3.4 formul\ae. Last but not least, this entails that Lemma 3.5 is Universal, too. 

\vfill
\eject

\par				
\noindent{\bf 3.2. Wintner Orthogonal Decompositions: $(P)$ and $)P($ properties} 
\bigskip
\par
\noindent
\par
\noindent
We need to bound series over $(P)$ and, also, to analyze counting functions both in $(P)$ and in $)P($. 
\par
Next Lemma provides these informations, coming from Lemma 2 and Lemma 3 in [C2].
\smallskip 
\par
\noindent{\bf Lemma 3.6}.
%
%
 {\stampatello (Converging series in $(P)$  \&  counting $P-$smooth and $P-$sifted numbers)}
\par
\noindent
{\it Let } $x\ge 1$ {\it be a} {\stampatello real number,} $P\in \Primes$ {\it and } $\delta\in \R, \delta>0$ {\it be } {\stampatello fixed}. {\it Then}
$$
\sum_{n\in (P)}n^{-\delta}\ll_{\delta,P} 1, 
\leqno{(1P)}
$$
$$
\sum_{{n\in (P)}\atop {n\le x}}1\ll_{\varepsilon,P} x^{\varepsilon},
\quad
\forall \varepsilon>0 ,
\leqno{(2P)}
$$
$$
\sum_{{n\in )P(}\atop {n\le x}}1=\prod_{p\le P}\left(1-{1\over p}\right)x+O_P(1). 
\leqno{(3P)}
$$
\smallskip
\par
\noindent{\bf Proof}. Any $n\in (P)$ can be {\stampatello represented as} \enspace $n=p_{1}^{K_1}\cdots p_{r}^{K_r}$, where \enspace $2=p_1<p_2<\cdots <p_r=P$ are {\stampatello consecutive prime numbers} and $K_j\in \N_0$, $\forall j\le r$; hence, $(1P)$ follows from: 
$$
\sum_{n\in (P)}n^{-\delta}=\sum_{K_1\ge 0}\left(p_1^{-\delta}\right)^{K_1}\cdots \sum_{K_r\ge 0}\left(p_r^{-\delta}\right)^{K_r}
={1\over 1-p_1^{-\delta}}\cdots {1\over 1-p_r^{-\delta}}
=\prod_{p\le P}{1\over 1-p^{-\delta}}
\ll_{\delta,P} 1. 
$$
\hfill \QED
\par
\noindent
The same representation gives $(2P)$ : $\forall \varepsilon>0$, from above $(1P)$, 
$$
\sum_{{n\in (P)}\atop {n\le x}}1\le \sum_{{n\in (P)}\atop {n\le x}}{{x^{\varepsilon}}\over {n^{\varepsilon}}}
\le x^{\varepsilon} \sum_{n\in (P)}n^{-\varepsilon}
\ll_{\varepsilon,P} x^{\varepsilon}. 
$$
\hfill \QED
\par
\noindent
Saying \enspace $n\in )P($\enspace is equivalent to : $(n,Q(P))=1$, where we abbreviate here $Q(P)={\displaystyle \prod_{p\le P}p }$; this coprimality condition is detected by the sum (see Eratosthenes-Legendre Sieve [HaRi]) \enspace ${\displaystyle \sum_{{d|n}\atop {d|Q(P)}}\mu(d) }$ : exchanging sums, 
$$
\sum_{{n\in )P(}\atop {n\le x}}1=\sum_{d|Q(P)}\mu(d)\sum_{{n\le x}\atop {n\equiv 0(\!\!\bmod d)}}1
=\sum_{d|Q(P)}\mu(d)\left({x\over d}+O(1)\right)
=x\sum_{d|Q(P)}{\mu(d)\over d}+O\left(\sum_{d|Q(P)}\mu^2(d)\right)
$$
\par
\noindent
and recalling {\stampatello Euler product}
$$
\sum_{d|n}{\mu(d)\over d}=\prod_{p|n}\left(1-{1\over p}\right)
\enspace \Longrightarrow \enspace 
\sum_{d|Q(P)}{\mu(d)\over d}=\prod_{p\le P}\left(1-{1\over p}\right)
$$
\par
\noindent
together with 
$$
\sum_{d|Q(P)}\mu^2(d)=\prod_{p\le P}2
=2^{\pi(P)}=O_P(1),
$$
\par
\noindent
we get $(3P)$.\hfill \square 
\medskip
\par
\noindent{\bf Remark 4}.
%
%
 For more details, in previous Lemma Proof ingredients, see $\S1$ in [C5] : Facts.\hfill $\diamond$

\vfill
\eject

\par				
\noindent{\bf 3.3. When are Carmichael's coefficients also Wintner's coefficients?} 
\bigskip
\par
\noindent
The answer to this question comes from a property that we expose here, in next result, and can be proved exactly like next Lemma 3.7, regarding the $P-$smooth version for Carmichael's \& Wintner's \lq \lq partial sums\rq \rq: 
\smallskip
\par
\noindent{\bf Proposition 3.2}.
%
%
 {\stampatello (Difference of Carmichael's and Wintner's partial sums)}
\par
\noindent
{\it Let } $F:\N \rightarrow \C$. {\it Then, } $\forall q\in \N$, $\forall x\in \R$, $x\ge 1$, 
$$
{1\over {\varphi(q)x}}\sum_{a\le x}c_q(a)F(a)=\sum_{{d\le x}\atop {d\equiv 0(\!\!\bmod q)}}{{F'(d)}\over d}
 +O_q\left({1\over x}\sum_{d\le x}|F'(d)|\right).
$$
\par
\noindent
The Proof mimics Lemma 3.7's one. 
\par
We saw above, before stating Wintner's Dream Theorem, Corollary 2.3, that $\CarT\, F = \WinT\, F$ follows when $F$ satisfies $\WA$ and this is contained in the quoted Wintner's Criterion. However, in [C3], we prove this, by equation $(5)$ there : actually, above formula for partial sums; then, soon after quoted $(5)$, we prove in [C3] that $\WA$ implies the condition: 
$$
\lim_x {1\over x}\sum_{d\le x}|F'(d)| = 0,
\leqno{\ETD}
$$
\par
\noindent
name abbreviating \lq \lq {\stampatello Eratosthenes Transform Decay}\rq \rq; namely, the {\bf Classic Mean Value} (compare [C3] again), of Eratosthenes Transform's modulus: $|F'|$, vanishes. See that, even if $\WA$ implies $\ETD$, this doesn't imply $\WA$ (a simple counterexample's $F$ with: $F'(1)=1$ and $F'(d)=1/\log(d)$, $\forall d>1$). 
\par
Summarizing, $\WA$ implies $\CarT\,F = \WinT\, F$; but this should hold under the condition $\ETD$ above, that is strictly weaker than $\WA$. However, while $\WA$, of course, gives $\exists \WinT\, F$, with $\ETD$ we need to know, furthermore, that $\exists \WinT\, F$ OR that $\exists \CarT\, F$, to conclude $\CarT\,F = \WinT\, F$. 
\par
We come, now, to the version for $P-$smooth counterparts of Carmichael's \& Wintner's coefficients. 
\par
\noindent
Next result tightly links Carmichael and Wintner $P-$smooth coefficients together, for a general arithmetic function $F$, once fixed $P\in \Primes$, estimating the remainder of any $q-$th $P-$smooth difference of such coefficients, for a kind of their $x-$partial sums. 
\par
Next Lemma Proof comes from exponential sums elementary bounds. 
\smallskip
\par
\noindent{\bf Lemma 3.7}.
%
%
 {\stampatello (Difference of $P-$smooth Carmichael's and $P-$smooth Wintner's partial sums)}
\par
\noindent
{\it Let } $F:\N \rightarrow \C$ {\stampatello and choose } $P\in \Primes$. {\it Then, } $\forall q\in \N$, $\forall x\in \R$, $x\ge 1$, 
$$
{1\over {\varphi(q)x}}\sum_{a\le x}c_q(a)\sum_{{d|a}\atop {d\in (P)}}F'(d)=\sum_{{d\le x}\atop {{d\in (P)}\atop {d\equiv 0(\!\!\bmod q)}}}{{F'(d)}\over d}
 +O_q\left({1\over x}\sum_{{d\le x}\atop {d\in (P)}}|F'(d)|\right).
$$
\smallskip
\par
\noindent{\bf Proof}. The LHS above is, exchanging sums, 
$$
{1\over {\varphi(q)x}}\sum_{{d\le x}\atop {d\in (P)}}F'(d)\sum_{m\le {x\over d}}c_q(dm)={1\over x}\sum_{{d\le x}\atop {{d\in (P)}\atop {d\equiv 0(\!\!\bmod q)}}}F'(d)\left({x\over d}+O(1)\right)
 +O\left(\sum_{{d\le x}\atop {{d\in (P)}\atop {d\not\equiv 0(\!\!\bmod q)}}}{{|F'(d)|}\over {\varphi(q)x}}\sum_{j\in \Z_q^*}\left|\sum_{m\le {x\over d}}e_q(jdm)\right|\right)
$$
\par
\noindent
which, thanks to [D,ch.26] that implies: 
$$
d\not \equiv 0(\bmod q)
\enspace \Longrightarrow \enspace 
\sum_{m\le {x\over d}}e_q(jdm)\ll \left\Vert {{jd}\over q}\right\Vert^{-1}
 = {1\over {\left\Vert {{j\cdot d/(d,q)}\over {q/(d,q)}}\right\Vert}}, 
$$
\par
\noindent
{\stampatello recalling} \enspace $\Vert \alpha\Vert \defineq \min_{n\in \Z}|\alpha-n|$, $\forall \alpha \in \R$, gives : LHS equals, abbreviating \enspace $q':=q/(d,q)\le q/2$, \enspace {\stampatello and changing variable} from $j\in \Z_q^*$ to $j'\in \Z_{q'}^*$ , 
$$
{1\over {\varphi(q)x}}\sum_{{d\le x}\atop {d\in (P)}}F'(d)\sum_{m\le {x\over d}}c_q(dm)
=\sum_{{d\le x}\atop {{d\in (P)}\atop {d\equiv 0(\!\!\bmod q)}}}{{F'(d)}\over d}
 +O\left({1\over x}\sum_{{d\le x}\atop {d\in (P)}}|F'(d)|\right)
  +O\left({1\over x}\sum_{{d\le x}\atop {d\in (P)}}{{|F'(d)|}\over {\varphi(q)}}\cdot{q\over {q'}}\cdot\sum_{j'\le {{q'}\over 2}}{{q'}\over {j'}}\right)
$$
\par
\noindent
{\stampatello in all}, with above remainder.\hfill $\square$ 

\vfill
\eject

\par				
\noindent{\bf 4. Odds \& Ends} 
\bigskip
\par
\noindent

\medskip

\par
\noindent{\bf 4.1 Smooth-twisted {\stampatello Orthogonality} of {\stampatello Ramanujan Sums}} 
\bigskip
\par
\noindent
Next result is a kind of NEW orthogonality property of Ramanujan sums, that we found in [C2] (3rd version!), compare Proposition 2 there. We adapt it here, to present applications. 
\par
Next Lemma gives a more suitable expression of this quoted result. 
\smallskip
\par
\noindent{\bf Lemma 4.1}.
%
%
 {\stampatello (Smooth$-$Twisted Orthogonality of Ramanujan Sums)}
\par
\noindent
{\it Let } $q,\ell\in (P)$, {\stampatello for a} {\stampatello fixed} $P\in \Primes$. {\it Then}
$$
{1\over {\varphi(\ell)
{\displaystyle \prod_{p\le P} }
\left(1-{1\over p}\right)^{-1}}}
 \sum_{t\in (P)}{{c_{\ell}(t)c_q(t)}\over t}
=
\1_{q=\ell}\thinspace . 
$$
\smallskip
\par
\noindent{\bf Proof}. We only need [C2,Proposition2], after showing: 
$$
\sum_{t\in (P)}{1\over t}=\prod_{p\le P}\sum_{K=0}^{\infty}{1\over {p^K}}
 =\prod_{p\le P}{1\over {1-{1\over p}}}
  =\prod_{p\le P}\left(1-{1\over p}\right)^{-1}.
$$ 
\par
\noindent
\enspace \hfill $\square$ 

\medskip

\par
\noindent{\bf 4.2 Square-free supported $F$ and their $F'$} 
\bigskip
\par
\noindent
We need to know why our expectation that s-f-supported a.f.s should have a cube-free supported Eratosthenes Transform is confirmed, so we do it in next result. In which we will use: $d$ cube-free $\Leftrightarrow $ $d/\kappa(d)$ square-free. 
\par
All of this is in next, very general and useful {\stampatello Lemma}. 
\smallskip
\par
\noindent{\bf Lemma 4.2}.
%
%
 {\stampatello (Eratosthenes Transform of square-free supported arithmetic functions)}
\par
\noindent
{\it Let } $F:\N \rightarrow \C$ {\it be} {\stampatello any arithmetic function}. {\it Then}
$$
\left(\mu^2 F\right)'(d)=\mu(d)\sum_{t|d}\mu(t)F(t)
+(1-\mu^2(d))\mu\left({d\over {\kappa(d)}}\right)\sum_
{
{t|\kappa(d)}\atop 
{
\left(
{{\kappa(d)}\over t},{d\over {\kappa(d)}}
\right)
=1
}
}
\mu\left({{\kappa(d)}\over t}\right)F(t),
\quad
\forall d\in \N.
$$
\par
\noindent{\bf Remark 5}.
%
%
 In the RHS, first term is on {\stampatello square-free} $d$, while second term's on $d$ {\stampatello not s-f, but cube-free}. In all, this proves, in particular, that {\stampatello any s-f-supported $F$ has cube-free-supp.$^{\hbox{\stampatello ed}}$} $F'$.\hfill $\diamond$
\smallskip
\par
\noindent
\par
\noindent{\bf Proof}. Use {\stampatello Eratosthenes Transform definition \& the property}: $d$ {\stampatello s-f divisors are} $\kappa(d)$ {\stampatello divisors} 
$$
\left(\mu^2 F\right)'(d)=\sum_{t|d}\mu^2(t)F(t)\mu\left({d\over t}\right)
 =\sum_{t|\kappa(d)}F(t)\mu\left({d\over {\kappa(d)}}\cdot {{\kappa(d)}\over t}\right)
  =\mu\left({d\over {\kappa(d)}}\right)\sum_{{t|\kappa(d)}\atop {(\kappa(d)/t,d/\kappa(d))=1}}\mu\left({{\kappa(d)}\over t}\right)F(t);
$$
\par
\noindent
distinguishing {\stampatello two incompatible cases} $\mu^2(d)=1$ {\stampatello and} $1-\mu^2(d)=1$, we get 
$$
\left(\mu^2 F\right)'(d)=\mu^2(d)\cdot \mu(d)\sum_{t|d}\mu(t)F(t)
+(1-\mu^2(d))\mu\left({d\over {\kappa(d)}}\right)
 \sum_{{t|\kappa(d)}\atop {\left({{\kappa(d)}\over t},{d\over {\kappa(d)}}\right)=1}}\mu\left({{\kappa(d)}\over t}\right)F(t),
$$
\par
\noindent
because \enspace $\mu^2(d)=1$ \enspace means $d$ square-free, whence \enspace $\mu(\kappa(d)/t)=\mu(d/t)=\mu(d)/\mu(t)=\mu(d)\mu(t)$, $\forall t|d$.\hfill $\square$ 
\medskip
\par
Before next section, we profit here to say that in our paper [CG2] we give a glance to the interactions of analytic and arithmetic aspects, for classic Ramanujan expansions. Needless to say, we would like to do the same (for example a classification of multiplicative coefficients), for the Smooth Ramanujan expansions. Like our first example $G=\1$, in $\S1$, shows, the coefficients for the two summation methods might be dramatically different, in fact. Last but not least : are there other interesting summation methods beyond the classic and our smooth one? 

\vfill
\eject

\par				
\noindent{\bf 5. Appendix} 
\medskip
We are concluding our short tour about elementary results for general arithmetic functions, regarding Smooth Ramanujan Expansions, both local and global: with a generalization for the Wintner Orthogonal Decomposition; Theorem 2.1 above works for the Wintner coefficients : these are, in turn, what Wintner Himself called \lq \lq Eratosthenian Averages\rq \rq: His Book's [W] name ! Then, we may generalize the Proof of Theorem 2.1 to a kind of more general \lq \lq Wintner Average\rq \rq, not only for $F'$, the Eratosthenes Transform. This we do in next subsection.  

\bigskip

\par
\noindent{\bf 5.1 Wintner Averages} 
\medskip
\par
\noindent
Given an arithmetic function, written like a sequence of complex numbers, 
$$
\{ a_n\}_{n\in \N}\thinspace , 
$$
\par
\noindent
whenever, of course, following series converge pointwise in all $q\in \N$,
$$
A_q \defineq \sum_{n\equiv 0(\!\!\bmod q)}{{a_n}\over n}\thinspace ,
$$
\par
\noindent
we call: $\{ A_q\}_{q\in \N}$ {\stampatello the Wintner Average of} our arithmetic function $\{ a_n\}_{n\in \N}$. 
\medskip
\par
\noindent{\bf Remark 6}.
%
%
 It may seem strange, at first sight, not to see $F'$, in these averages: it may seem a kind of considering, say, $A_q=\Win_q\,(F\ast 1)$, with $\ast$ Dirichlet product [T]. However, compare next Remark.\hfill $\diamond$
\medskip
\par
Following holds.
\smallskip
\par
\noindent{\bf Proposition 3.3}.
%
%
 {\stampatello (Wintner Orthogonal Decomposition with Wintner Averages)}
\par
\noindent
{\it Let } $\{ a_n\}_{n\in \N}$ {\it have the} {\stampatello Wintner Average}, {\it say}
$$
A : q\in \N \mapsto A_q \EqByDef \sum_{n\equiv 0(\!\!\bmod q)}{{a_n}\over n}\in \C, 
$$
{\it we defined above. Fix a prime} $P$. {\it Then}
$$
\forall d\in \N, 
\qquad
\sum_{{r\in )P(}\atop {r>1}}{{a_{dr}}\over r}\defineq \lim_x \sum_{{r\in )P(}\atop {1<r\le x}}{{a_{dr}}\over r}\in \C
$$
\par
\noindent
{\it is a series, converging pointwise in } $d$, $\forall d\in \N$, {\it called the } {\stampatello Irregularity} {\it of our } $\{ a_n\}_{n\in \N}$, {\stampatello with argument } $d$, {\stampatello over the prime } $P$. {\it Varying $d\in \N$, we get the arithmetic function}
$$
{\rm Irreg}^{(P)}\,a : d\in \N \mapsto \sum_{{r\in )P(}\atop {r>1}}{{a_{dr}}\over r}\in \C. 
$$
\par
\noindent
{\it In all, } $\exists \{ A_q\}_{q\in \N}$ {\it entails } $\exists {\rm Irreg}^{(P)}\,a$, $\forall P\in \Primes$; {\it and} {\stampatello Wintner Average Decomposition of } $\{ a_n\}_{n\in \N}$ : $\forall P\in \Primes$,
$$
a_d=d\sum_{K\in (P)}\mu(K)A_{dK} - \sum_{{r\in )P(}\atop {r>1}}{{a_{dr}}\over r},
\qquad
\forall d\in \N. 
\leqno{\WAD_{a,P}}
$$
\smallskip
\par
\noindent{\stampatello Proof} works with $a_d$ instead of $F'(d)$, $\forall d\in \N$, because we mimic Theorem 2.1 Proof.\hfill $\square$
\medskip
\par
\noindent{\bf Remark 7}.
%
%
 The {\stampatello Irregular Series for} $F$ is {\stampatello not the Irregularity of} our $F$, {\stampatello but that of} our $F'$.\hfill $\diamond$

\vfill
\eject

\par				
\centerline{\stampatello Bibliography}

\bigskip

\item{[C1]} G. Coppola, {\sl An elementary property of correlations}, Hardy-Ramanujan J. {\bf 41} (2018), 65--76.
\smallskip
\item{[C2]} G. Coppola, {\sl A smooth shift approach for a Ramanujan expansion}, arXiv:1901.01584v3 (3rd version)
\smallskip
\item{[C3]} G. Coppola, {\sl Recent results on Ramanujan expansions with applications to correlations}, Rend. Sem. Mat. Univ. Pol. Torino {\bf 78.1} (2020), 57--82. 
\smallskip
\item{[C4]} G. Coppola, {\sl A smooth summation of Ramanujan expansions}, arXiv:2012.11231v8 (8th version)
\smallskip
\item{[C5]} G. Coppola, {\sl General elementary methods meeting elementary properties \thinspace of \thinspace correlations}, {\tt available}\break{\tt online at} \enspace arXiv:2309.17101 (2nd version)
\smallskip
\item{[CG1]} G. Coppola and L. Ghidelli, {\sl Multiplicative Ramanujan coefficients of null-function}, arXiv:2005.14666v2 (2nd Version) 
\smallskip
\item{[CG2]} G. Coppola and L. Ghidelli, {\sl Convergence of Ramanujan expansions, I [Multiplicativity on Ramanujan clouds]}, arXiv:2009.14121v1 
\smallskip
\item{[CM]} G. Coppola and M. Ram Murty, {\sl Finite Ramanujan expansions and shifted convolution sums of arithmetical functions, II}, J. Number Theory {\bf 185} (2018), 16--47. 
\smallskip
\item{[D]} H. Davenport, {\sl Multiplicative Number Theory}, 3rd ed., GTM 74, Springer, New York, 2000. 
\smallskip
\item{[De]} H. Delange, {\sl On Ramanujan expansions of certain arithmetical functions}, Acta Arith., {\bf 31} (1976), 259--270.
\smallskip
\item{[HaRi]} H. Halberstam and H.E. Richert, {\sl Sieve Methods}, London Mathematical Society Monographs, No. 4. {\it Academic Press, London-New York}, 1974. xiv + 364 pp.
\smallskip
\item{[H]} G.H. Hardy, {\sl Note on Ramanujan's trigonometrical function $c_q(n)$ and certain series of arithmetical functions}, Proc. Cambridge Phil. Soc. {\bf 20} (1921), 263--271.
\smallskip
\item{[IKo]} H. Iwaniec and E. Kowalski, {\sl Analytic Number Theory}, American Mathematical Society Colloquium Publications, 53. American Mathematical Society, Providence, RI, 2004. xii+615pp. \hfil ISBN:0-8218-3633-1
\smallskip
\item{[K]} J.C. Kluyver, {\sl Some formulae concerning the integers less than $n$ and prime to $n$}, Proceedings of the Royal Netherlands Academy of Arts and Sciences (KNAW), {\bf 9(1)} (1906), 408--414. 
\smallskip
\item{[M]} M. Ram Murty, {\sl Ramanujan series for arithmetical functions}, Hardy-Ramanujan J. {\bf 36} (2013), 21--33. Available online 
\smallskip
\item{[R]} S. Ramanujan, {\sl On certain trigonometrical sums and their application to the theory of numbers}, Transactions Cambr. Phil. Soc. {\bf 22} (1918), 259--276.
\smallskip
\item{[ScSp]} W. Schwarz and J. Spilker, {\sl Arithmetical Functions}, Cambridge University Press, 1994.
\smallskip
\item{[T]} G. Tenenbaum, {\sl Introduction to Analytic and Probabilistic Number Theory}, Cambridge Studies in Advanced Mathematics, {46}, Cambridge University Press, 1995. 
\smallskip
\item{[W]} A. Wintner, {\sl Eratosthenian averages}, Waverly Press, Baltimore, MD, 1943. 

\bigskip
\bigskip
\bigskip

\par
\leftline{\tt Giovanni Coppola - Universit\`{a} degli Studi di Salerno (affiliation)}
\leftline{\tt Home address : Via Partenio 12 - 83100, Avellino (AV) - ITALY}
\leftline{\tt e-mail : giocop70@gmail.com}
\leftline{\tt e-page : www.giovannicoppola.name}
\leftline{\tt e-site : www.researchgate.net}

\bye